\documentclass[12pt]{article}
\usepackage[a4paper]{geometry}
\usepackage{amsthm,amsmath,amssymb,enumerate,color,srcltx}
\usepackage{natbib}

\usepackage{graphicx}	% insert graph
\usepackage{float}

\usepackage{booktabs}   % insert table
\usepackage{array}
\usepackage[english]{babel}
\usepackage{lscape}

\usepackage[colorlinks]{hyperref}

\newtheorem{theorem}{Theorem}[section]

\newtheorem{assumption}{Assumption}[section]
\newtheorem{lemma}{Lemma}[section]
\newtheorem{proposition}{Proposition}[section]

\theoremstyle{remark}

\newcommand{\Rset}{{\mathbb{R}}} %%% real numbers
\newcommand{\Zset}{{\mathbb{Z}}} %%% integer numbers

\newcommand{\esp}{{\mathbb E}}
\newcommand{\pr}{{\mathbb P}}
\newcommand{\var}{{\mathrm{var}}}
\newcommand{\cov}{{\mathrm{cov}}}
\newcommand{\rme}{{\mathrm{e}}}
\newcommand{\rmi}{{\mathrm{i}}}
\newcommand{\rmd}{{\mathrm{d}}}

\newcommand{\plim}{{\ \stackrel{\mathbb{P}}\to} \ }

\newcommand{\dlim}{{\to}}
\newcommand{\func}{\ \Rightarrow \ }
\newcommand{\fidi}{\ \stackrel{fi.di.}\longrightarrow\ }
\newcommand{\iid}{\ \stackrel{i.i.d.}\sim\ }
\newcommand{\lip}{\mathrm{Lip}}
\newcommand{\mcg}{\mathcal{G}}

\newcommand{\newstat}{T}
\newcommand{\limproc}{A}
\newcommand{\reste}{{\rm rest}}

\numberwithin{equation}{section}

\begin{document}

\title{Drift in Transaction-Level Asset Price Models}

\author{Wen Cao\footnote{Stern School of Business, New York University, Henry
    Kaufman Management Center, 44 West Fourth Street, New York, NY 10012, USA,
    email: {\tt wcao@stern.nyu.edu}} \and Clifford Hurvich\footnote{Stern School
    of Business, New York University, Henry Kaufman Management Center, 44 West
    Fourth Street, New York, NY 10012, USA, email: {\tt churvich@stern.nyu.edu}}
  \and Philippe Soulier\footnote{Universit\'e Paris Ouest, 200 avenue de la
    R\'epublique, 92001 Nanterre cedex, France, email: {\tt
      philippe.soulier@u-paris10.fr. The research of Philippe Soulier was
      partially supported by the grant ANR-08-BLAN-0314-02. } \newline The
    authors thank participants of the 2012 NBER-NSF Time Series Conference for
    helpful comments.} }

\maketitle

\begin{abstract}
  We study the effect of drift in pure-jump transaction-level models for asset prices in continuous
  time, driven by point processes. The drift is assumed to arise from a nonzero mean in the
  efficient shock series. It follows that the drift is proportional to the driving point process
  itself, i.e. the cumulative number of transactions. This link reveals a mechanism by which
  properties of intertrade durations (such as heavy tails and long memory) can have a strong impact
  on properties of average returns, thereby potentially making it extremely difficult to determine
  long-term growth rates or to reliably detect an equity premium.  We focus on a basic univariate
  model for log price, coupled with general assumptions on the point process that are satisfied by
  several existing flexible models, allowing for both long memory and heavy tails in
  durations. Under our pure-jump model, we obtain the limiting distribution for the suitably
  normalized log price. This limiting distribution need not be Gaussian, and may have either finite
  variance or infinite variance.  We show that the drift can affect not only the limiting
  distribution for the normalized log price, but also the rate in the corresponding normalization.
  Therefore, the drift (or equivalently, the properties of durations) affects the rate of
  convergence of estimators of the growth rate, and can invalidate standard hypothesis tests for
  that growth rate.  As a remedy to these problems, we propose a new ratio statistic which behaves
  more robustly, and employ subsampling methods to carry out inference for the growth rate.  Our
  analysis also sheds some new light on two longstanding debates as to whether stock returns have
  long memory or infinite variance.
\end{abstract}

\section{Introduction}

In recent years, transaction-level data on financial markets has become increasingly available, and
is now often used to make trading decisions in real time. Such data typically consist of the times
at which transactions occurred, together with the price at which the transaction was executed, and
may include other concomitant variables ("marks") such as the number of shares traded. Our focus
here is on actual transactions rather than quotes, but regardless of which type of event is being
considered it is important to recognize that a useful framework for modeling and analyzing such data
is that of marked point processes rather than, say, time series in discrete time. Though time series
are typically provided for further analysis, such as daily (or high frequency) stock returns, these
inevitably involve aggregation and entail a loss of information that may be crucial for trading and
perhaps even for risk management and portfolio selection.

The perspective of asset prices as (marked) point processes has a long history in the financial and
econometric literature. For example, \citet{scholes:williams:1977} allowed for a compound Poisson
process. However, such a model is at odds with the stylized fact that time series of financial
returns exhibit persistence in volatility.  Recent interest in the point process approach to
modeling transaction-level data was spurred by the seminal paper of \citet{engle:russell:1998}, who
proposed a model for inter-trade durations.  Other work on modeling transaction-level data as point
processes and/or constructing duration models includes that of \citet{prigent:2001},
\citet{bowsher:2007}, \citet{billingsley:1968}, \citet{hautsch:2012},
\citet{bacry:delattre:hoffmann:muzy::2011}, \citet{deo:hurvich:soulier:wang:2009},
\citet{deo:hsieh:hurvich:2010}, \citet{hurvich:wang:2009},
\citet{aue:hurvich:horvath:soulier:2014}, \citet{shenai:2012},
\citet{chen:diebold:schorfheide}.

Nevertheless, it must be recognized that time series of asset returns in
discrete (say, equally-spaced) time are still in widespread use, and indeed may
be the only recorded form of the data that encompasses many decades. Such long
historical series are of importance for understanding long-term trends (a prime
focus of this paper) and, arguably, for a realistic assessment of risk. So given
the ubiquitous nature of the time series data but also keeping in mind the
underlying price-generating process that occurred at the level of individual
transactions, it is important to make sure that transaction-level models obey
the stylized facts, not only for the intertrade durations but also for the
lower-frequency time series.

It has been observed empirically that time series of financial returns are
weakly autocorrelated (though perhaps not completely uncorrelated), while
squared returns or other proxies for volatility show strong autocorrelations
that decay very slowly with increasing lag, possibly suggesting long memory (see
\citet{andersen:bollerslev:diebold:2001}). It is also generally accepted that such time series show
asymmetries, such as a correlation between the current return and the next
period's squared return, and this effect (often referred to traditionally as the
"leverage effect") is addressed for example by the EGARCH model of
\citet{nelson:1991}. The average return often differs significantly from zero based
on a traditional $t$-test, possibly suggesting a linear trend in the series of
log prices. Meanwhile, \citet{deo:hsieh:hurvich:2010} found that intertrade
durations have long memory (this was also found by
\citet{chen:diebold:schorfheide}), and they investigated the possibility that the
durations have heavy tails.

One more fact that we wish to stress is that in continuous time, realizations of
series of transaction-based asset prices are step functions, since the price is
constant unless a transaction occurs.  Thus, we choose to focus on pure-jump
models for the log price (viewed as a time series in continuous time), driven by
a point process that counts the cumulative number of transactions. This is
equivalent to a marked point process approach where the points correspond to the
transaction times and the marks are the transaction-level return shocks, which
can have both an efficient component and a microstructure component.

Within this context, we will in this paper investigate the effect of drift
(modeled at the transaction level) on the behavior of very-long-horizon returns,
or equivalently, on the asymptotic behavior of the log price as time
increases. The drift is assumed to arise from a nonzero mean in the efficient
shock series. It follows that the drift is proportional to the driving point
process itself, i.e. the cumulative number of transactions. This link reveals a
mechanism by which properties of intertrade durations (such as heavy tails and
long memory) can have a strong impact on properties of average returns, thereby
potentially making it extremely difficult to determine long-term growth rates or
to reliably detect an equity premium.

We focus on a basic univariate model for log price, coupled with general
assumptions on the point process that are satisfied by several existing flexible
models, allowing for both long memory and heavy tails in durations. Under our
pure-jump model (which can capture all the stylized facts described above), we
obtain the limiting distribution for the suitably normalized log price. This
limiting distribution need not be Gaussian, and may have either finite variance
or infinite variance. The diversity of limiting distributions here may be
considered surprising, since our assumptions imply that the return shocks obey
an ordinary central limit theorem under aggregation across transactions ($i.e.$
in $transaction$ $time$ but generally not in $calendar$ $time$). We show that
the drift can affect not only the limiting distribution for the normalized log
price, but also the rate in the corresponding normalization. Therefore, the
drift (or equivalently, the properties of durations) affects the rate of
convergence of estimators of the growth rate, and can invalidate standard
hypothesis tests for that growth rate. Our analysis also sheds some new light on
two longstanding debates as to whether stock returns have long memory or
infinite variance.

The remainder of this paper is organized as follows. In Section \ref{sec:model},
we provide a simple univariate model for the log price, discuss the trend term
and present Proposition \ref{prop:weakconv-returns} on the limiting behavior of
the log price process, as determined by the properties of the point process. In
Section \ref{sec:Inference} we begin a study of statistical inference for the
trend, and obtain the behavior of the ordinary $t$-statistic under the null
hypothesis. We also propose a ratio statistic which behaves robustly, in that
(unlike the $t$-statistic) it converges in distribution under a broad range of
conditions. This motivates our discussion in Section \ref{sec:Subsampling} of
the use of subsampling based on the ratio statistic to conduct statistical
inference on the trend. We then present, in Section \ref{sec:examples}, a series
of examples of specific duration and point process models that have been
proposed in the literature, including the Autoregressive conditional duration (ACD) model of \citet{engle:russell:1998} and the
 long memory stochastic duration (LMSD) model of \citet{deo:hsieh:hurvich:2010}. These examples provide for great
diversity of the asymptotic distributions of sums of durations and therefore (by
Proposition \ref{prop:weakconv-returns} below) for the asymptotic distribution
of the log price. Simulations are presented in Section \ref{sec:simulation}. There, we consider the size and power for both the $t$-test and the test based on the new ratio statistic under three duration
models: $iid$ exponential, ACD and LMSD. The parameters of simulations are calibrated using empirical data. Section \ref{sec:data} provides a data analysis in which we employ the $t$-test and the test based on the new ratio statistic to gauge the strength of the equity premium. Section \ref{sec:discussion} provides a concluding discussion
on how our results may help to reconcile some longstanding debates. Proofs of
the mathematical results are provided in Section \ref{sec:math}.

\section{A Univariate Model for Log Price}
\label{sec:model} We start with a basic univariate pure-jump model for a log
price series $y$. Let $N_1$, $N_2$ be mutually independent point processes on
the real line. Points of $N_1+N_2$ correspond to transactions on the
asset\footnote {From a modeling perspective it may be desirable to instead have
  the points of $N_1+N_2$ correspond to other relevant trading events, such as
  "every fourth transaction", "a transaction that moves the price", etc. For
  simplicity and definiteness in the paper, we simply let $(N_1+N_2)(t)$ count
  actual transactions, but our theoretical results do not depend on this
  particular choice of the definition of an event.} and $(N_1+N_2)(t)$ is the
number of transactions in $(0,t]$. We set $y(0)=0$ and define, for $t \geq 0$,
\begin{align}
  \label{eq:UnivariateModel}
  y(t) = \sum_{k=1}^{N_1(t)} (\tilde e_{1,k}+\eta_{1,k}) + \sum_{k=1}^{N_2(t)}
  (\tilde e_{2,k}+\eta_{2,k})
\end{align}
where for $i=1,2$,
\[
\tilde e_{i,k} = \mu_i + e_{i,k}
\]
and the $\{e_{i,k}\}$, which are independent of $(N_1,N_2)$, are i.i.d.~with
zero mean and finite variance $\sigma_i^2$.  We assume that the microstructure
noise sequences $\{\eta_{i,k}\}$ are independent of the efficient shocks
$\{e_{1,k}\}$, $\{e_{2,k}\}$, but not necessarily of the counting process
$(N_1,N_2)$. This allows for a leverage effect. See
\citet{aue:hurvich:horvath:soulier:2014}.  Assumption \ref{hypo:micro-negligible}
below implies that the microstructure noise becomes negligible after
aggregation.

We assume that $\mu_1$ and $\mu_2$ are constants, not both zero. This model with
$\mu_1=\mu_2=0$, $\sigma_2^2=0$ and $\{\eta_{2,k}\}=0$ was considered in \citet
{deo:hurvich:soulier:wang:2009}, who showed that it can produce long memory in
the realized volatility. We generalize the model here to allow for two driving
processes, $N_1$ and $N_2$.  It follows from~(\ref{eq:UnivariateModel}) that
\begin{align}
  \label{eq:trend}
  y(t) = \mu_1 N_1(t) + \mu_2N_2(t) + \sum_{k=1}^{N_1(t)} (e_{1,k}+\eta_{1,k}) + \sum_{k=1}^{N_2(t)} (e_{2,k}+\eta_{2,k}) \; .
\end{align}
The quantity $\mu_1 N_1(t) + \mu_2N_2(t)$ can be viewed as a random drift term
in the log price $y(t)$. If, for example, $\mu_1 >0$ and $\mu_2 <0$ then we can
think of $N_1$ and $N_2$ as governing the activity of buyers and sellers,
respectively. An active period for $N_1$ puts upward pressure on the price,
while an active period for $N_2$ exerts downward pressure. One could envision
including other point processes to differentiate, say, between information
traders and noise traders, but we find that the model with just $N_1$ and $N_2$
is quite flexible and able to reproduce the stylized facts.

To simplify the discussion, we assume temporarily that $N_1$ and $N_2$ are
stationary point processes with intensities $\lambda_1$ and $\lambda_2$. Then
the expectation of the drift term is a linear trend, that is, $\esp [\mu_1
N_1(t) + \mu_2N_2(t)]=\mu^* t$, where
\[
\mu^* = \lambda_1 \mu_1 + \lambda_2 \mu_2 \,\,\, .
\]
Since we are modeling the log prices $y(t)$ as a pure-jump process the log price
is constant when no trading occurs. Unfortunately, the modified version of
(\ref{eq:trend}) in which the random drift term $\mu_1 N_1(t) + \mu_2N_2(t)$ is
replaced by the deterministic time trend $ct$ (where $c$ is a nonzero constant)
would not yield a pure-jump process. Nevertheless, it is quite reasonable from
an economic viewpoint to imagine that $\esp [y(t)]$ is a linear function of $t$,
to account for such phenomena as equity premia and inflation.  This is indeed
the case for Model~(\ref{eq:trend}) if it is assumed in addition that the
microstructure noise terms have zero mean under their Palm distribution (see
below), which implies $\esp [y(t)]=\mu^*t$. But in actual realizations of
$y(t)$, the trend is only impounded when a transaction occurs, due to the
nonzero means in~$\{\tilde e_{1,k}\}$ and/or $\{\tilde e_{2,k}\}$.

We now obtain conditions under which the two summations on the righthand side of
(\ref{eq:trend}) have mean zero, so that $\esp [y(t)]=\mu^*t$, still under the
assumption that $N_1$ and $N_2$ are stationary. Since the efficient shocks are
assumed independent of $N_1$ and $N_2$ and have mean zero, it follows that
$\esp\left[\sum_{k=1}^{N_i(t)} e_{i,k}\right] = 0$.  However, due to potential
leverage-type effects, even if one were to assume that $\esp[\eta_{i,k}]=0$,
this alone would not ensure that $\esp\left[\sum_{k=1}^{N_i(t)}
  \eta_{i,k}\right] = 0$. Therefore, if we want to interpret expectation of the
drift term as a linear trend in the log price, we must make the additional
assumption that $\esp\left[\sum_{k=1}^{N_i(t)} \eta_{i,k}\right] = 0$.  As we
will show, this assumption is equivalent to assuming that the
$E_i^0[\eta_{i,k}]=0$, where $P_i^0$ is the Palm probability measure associated
with the point process $N_i$, and $E_i^0$ is the expectation under the Palm
measure.

We briefly recall here (see \citet{deo:hurvich:soulier:wang:2009} for more
details) that in general there is no single measure under which both the
durations and the marked point process are stationary. If the durations for a
point process $N$ are stationary under a measure $P^0$, then we refer to $P^0$
as the Palm probability measure. This, if such a measure exists, would be the
appropriate measure for assessing properties of durations as a stationary
sequence. In such situations, there exists a corresponding measure $\pr$ under
which the point process is stationary. This would be the appropriate
distribution for assessing counts (the number of events occurring in equally
spaced intervals of time), and also the differences of $y(t)$ in
(\ref{eq:trend}), as a stationary time series in models that admit stationarity.

If $\{t_k,z_k,k\in\Zset\}$ are the (marked) points of a stationary (under
$\pr$) marked point process $N$ with finite intensity $\lambda$, then, by
\citet[Formula~1.2.9]{baccelli:bremaud:2003}, for all $t>0$
\begin{align*}
  \esp \left[ \sum_{k=1}^{N(t)} z_k \right] = \lambda t E^0[z_0] \; ,
\end{align*}
where $E^0$ is the expectation with respect to the Palm probability $P^0$.  If
the marks $\{z_k\}$ have zero mean under the Palm measure $P^0$, then $ \esp
\left[ \sum_{k=1}^{N(t)} z_k \right] =0$, even if, under $\pr$, it might happen
that $\esp[z_0] = \lambda E^0[t_1z_1] \ne 0$.

It follows that if the microstructure noise sequences $\{\eta_{i,k}\}$ have zero
mean under their respective Palm measures $P_i^0$ then in the model
(\ref{eq:trend}) under $\pr$ we have $\esp [y(t)] = \mu^*t$.

There are a variety of reasons why it may be considered unrealistic to assume
that the point process $(N_1,N_2)$ is stationary. Indeed, it is well-known that
the counts of financial transactions can show intraday seasonality (see, e.g.,
\citet{engle:russell:1998}, \citet{deo:hsieh:hurvich:2010}), and the transaction counts
will be guaranteed to be zero when the market is not open for trading.
Furthermore, regime changes are widely understood to be an important feature of
economic processes. We will consider a time deformation mechanism below which
allows for such effects. In the absence of stationarity we are no longer able to
compute $\esp[y(t)]$ in (\ref{eq:trend}), but we can still interpret $\mu^*$ as
a long term trend, i.e. $y(t)/t\plim \mu^*$, under the assumptions below. These
assumptions are stronger than needed for this purpose but will be used later to
derive limiting distributions for suitably rescaled versions of $y$.

In the sequel, we will consider a single probability measure $\pr$, without any requirement that
$(N_1,N_2)$ be stationary under $\pr$, unless specified otherwise.
% In the mathematical theory presented in this paper, all random variables and
% stochastic processes are defined on a single probability space
% $(\Omega,\mathcal{F},\pr)$. Expectation with respect to $\pr$ will be denoted by
% $\esp$ and $\var$ and $\cov$ will denote the variance and covariance with
% respect to $\pr$.
Convergence in $\pr$-probability will be denoted by~$\plim$, convergence in distribution under $\pr$
of sequences of random variables will be denoted by~$\dlim$ and convergence of finite dimensional
distributions of a sequence of stochastic processes will be denoted by~$\fidi$.
% We use $\func$ to denote weak convergence under $\pr$ in the space $\mathcal D([0,\infty))$ of
% left-limited, right-continuous (c\`adl\`ag) functions, endowed with Skorohod's $J_1$ topology.
% See \citet{billingsley:1968} or \citet{whitt:2002} for details about weak convergence in $\mathcal
% D([0,\infty))$. Whenever the limiting process is continuous, this topology can be replaced by the
% topology of uniform convergence on compact sets.

\begin{assumption}
  \label{hypo:weakconv-counts}
  The point processes $N_1$ and $N_2$ are mutually independent and there exists
  $\gamma \geq 1/2$ such that
  \begin{align}
    \label{eq:weakconv-counts}
    \{ n^{-\gamma} (N_i(nt)- \lambda_int),t\geq0\} \fidi \limproc_i
  \end{align}
  for $i=1,2$, as $n\rightarrow \infty$ where $\limproc_i$ are stochastic processes, at
  least one of which is nonzero.
\end{assumption}
We will give examples of point processes satisfying
Assumption~\ref{hypo:weakconv-counts} in Section~\ref{sec:examples}.  As will be
seen, the limiting processes $\limproc_i$ can be Gaussian or have infinite
variance, can have independent or dependent increments, and $\gamma$ can take
any value in $[1/2,1)$.  However, whatever $\limproc$ and $\gamma$, a consequence of
the convergence~(\ref{eq:weakconv-counts}) is that $N_i(t)/t\plim\lambda_i$ and
$t^{-1/2} \sum_{k=1}^{N_i(t)} e_{i,k}$ converges weakly to a Brownian motion.

Two distinct modeling approaches seem natural. One is to model the point process
directly as in \citet{bacry:delattre:hoffmann:muzy::2011} who use Hawkes
processes or \citet{delattre:robert:rosenbaum:2013} who use Cox processes.
Another approach consists of modeling the durations. This was done
by \citet{engle:russell:1998} who defined the ACD model and
\citet{deo:hsieh:hurvich:2010} who defined the LMSD model. These two approaches
are equivalent with respect to Assumption~\ref{hypo:weakconv-counts}, since
(\ref{eq:weakconv-counts}) with limit process $\limproc$ holds for a point process $N$
with intensity $\lambda$ and duration sequence $\{\tau_k\}$ if and only if
\begin{align}
  \label{eq:weakconv-durations}
  n^{-\gamma} \sum_{k=1}^{[n\cdot]} (\tau_k -\lambda^{-1}) \fidi
  -\lambda^{-(1+\gamma)} \limproc \; ,
\end{align}
where convergence holds in the sense of finite dimensional
distributions. See \citet[Theorem~7.3.1]{whitt:2002}.

\begin{assumption}
  \label{hypo:micro-negligible}
  The sequence of processes $n^{-1/2} \sum_{k=1}^{[n\cdot]} \eta_{i,k}$ converges in
  probability uniformly on compact sets to 0.
\end{assumption}

Assumptions~\ref{hypo:weakconv-counts} and~\ref{hypo:micro-negligible} imply that
$t^{-1/2}\sum_{k=1}^{N_i(t)} (e_{i,k}+\eta_{i,k})$ converges weakly to a Brownian motion and so the
microstructure terms are asymptotically negligible.

The following proposition, which follows easily from Assumptions \ref{hypo:weakconv-counts} and
\ref{hypo:micro-negligible}, gives the limiting distribution for suitably rescaled versions of the
log price process $y$. Recall that Assumption \ref{hypo:weakconv-counts} implies that
$N_i(t)/t\plim\lambda_i$. Let $\mu^*=\lambda_1\mu_1+\lambda_2\mu_2$.
\begin{proposition}
  \label{prop:weakconv-returns}
  Let Assumptions~\ref{hypo:weakconv-counts} and~\ref{hypo:micro-negligible} hold.
  Then $y(t)/t\plim\mu^*$. Moreover, as $n\rightarrow \infty$,
  \begin{align}
   \label{eq:conweak-gamma>1/2}
   \{n^{-\gamma} (y(nt) - \mu^* nt), t\geq0\} \fidi \mu_1 \limproc_1 + \mu_2 \limproc_2 \; ,
  \end{align}
  if $\gamma > 1/2$ and
  \begin{align}
   \label{eq:conweak-1/2}
   \{ n^{-1/2} (y(nt) - \mu^* nt) , t\geq0\} \fidi \mu_1 \limproc_1 + \mu_2 \limproc_2 + \sigma B
  \end{align}
  if $\gamma=1/2$ or if $\mu_1=\mu_2=0$, where $B$ is a standard Brownian motion independent of
  $\limproc_i$ and $\sigma$ is a positive constant.
\end{proposition}

Next, we discuss the use of time deformation. This discussion is taken from
\citet{aue:hurvich:horvath:soulier:2014}. Let $f$ be a deterministic or random function such that
$f$ is nondecreasing and has c\`adl\`ag paths with probability one.  Let $\tilde N$ be a point
process and define
\[
N(t)=\tilde N(f(t)) \; .
\]
We do not require that $\tilde N$ be stationary, but one useful application of time deformation is
to start with a stationary process $\tilde N$ and then deform it as above to obtain a nonstationary
process $N$. If the function $f$ is random, we assume moreover that it is independent of $\tilde
N$. The use of the time-varying intensity function $f$ may render the counting process $N$
nonstationary even if $\tilde N$ is a stationary point process. Since it is possible that $f$ has
(upward) jumps, the point process $N$ may also not be simple even though $\tilde N$ is simple.  We
now show, however, that if $\tilde N$ satisfies Assumption~\ref{hypo:weakconv-counts}, then so does
the time-deformed $N$ under some restrictions on $f$.

\begin{lemma}
  \label{lem:TimeDeformation-nonstationaire}
  Assume that $f$ is a nondecreasing (random) function such that $t^{-1}f(t)
  \plim \delta \in(0,\infty)$ and
\[
\sup_{t\geq 0} |f(t) - f(t^-)| \leq C
\]
with probability one, where $C \in (0,\infty)$ is a deterministic constant. Let $\tilde N$ be a
point process such that Assumption~\ref{hypo:weakconv-counts} holds for some $\tilde
\lambda\in(0,\infty)$ and $\gamma>0$. Let $N$ be the counting process defined by $N(\cdot)=\tilde
N(f(\cdot))$.  Then Assumption~\ref{hypo:weakconv-counts} holds for $N$ with
$\lambda=\tilde\lambda\delta$.
\end{lemma}

The function $f$ is used to speed up or slow down the trading clock.  To incorporate dynamic
intraday seasonality in volatility, the same time deformation can be used in each trading period (of
length, say, $T$), assuming that $f(t)$ has a periodic derivative (with period $T$ and with
probability one), for example, $f(t)=t+.5\sin(2\pi t/T)$.  Fixed nontrading intervals, say, $t \in
[T_1,T_2)$, could be accommodated by taking $f(t)=f(T_1)$ for $t \in [T_1,T_2)$ so that $f(t)$
remains constant for $t$ in this interval, and then taking $f(T_2) > f(T_1)$ so that $f(t)$ jumps
upward when trading resumes at time $T_2$. The jump allows for the possibility of one or more
transactions at time $T_2$, potentially reflecting information from other markets or assets that did
trade in the period $[T_1,T_2)$.  Since it is possible that $f$ has (upward) jumps, $N$ may not
be simple even though $\tilde N$ is simple.

\section{Statistical inference for the trend}
\label{sec:Inference}

For integer $k$, (assuming a time-spacing of 1 without loss of generality) we define the
calendar-time returns as $r_k = y(k)-y(k-1)$ and the average return over a time period that extends
from time $0$ to time $n$ as $\bar r_n = n^{-1} y(n)= n^{-1} \sum_{k=1}^n r_k$. We can think of $n$
as the length of the time interval spanned by the observations, or (since we have assumed a
time-spacing of 1) as the sample size, $i.e.$, the number of observed returns. It follows from
Proposition \ref{prop:weakconv-returns} that if $\gamma > 1/2$, $\bar r_n-\mu^*$ will not be
$O_p(n^{-1/2})$, making it difficult to accurately estimate growth rates based on the data set
$\{r_k\}_{k=1}^n$.

We recall that model (\ref{eq:trend}) under stationarity of $N_1$ and
$N_2$ implies that $\esp[y(n)] =
  \mu^* n$, so that the growth rate per unit time is $\mu^*$, but
  that under the more general assumptions \ref{hypo:weakconv-counts}
  and \ref{hypo:micro-negligible} $\mu^*$ may still be viewed as a long-term
  growth rate, since $y(n)/n\plim \mu^*$.

We consider the problem of statistical inference for $\mu^*$.  We focus on
testing a null hypothesis of the form $H_0 : \mu^* = \mu_0^*$ based on $\bar r_n$.

\subsection{The $t$ statistic}
\label{subsec:tstat}
The corresponding $t$-statistic for testing $H_0$ is
\[
t_n = n^{1/2}(\bar r_n - \mu_0^*)/s_n,
\]
where
\[
s_n^2=(n-1)^{-1} \sum_{j=1}^n (r_j - \bar r_n)^2.
\]

The sample variance $s_n^2$ consistently estimates a positive
constant, under suitable moments assumptions.

\begin{lemma}
   \label{lem:consistency-s}
   Under Assumptions~\ref{hypo:weakconv-counts} and~\ref{hypo:micro-negligible},
   if $N_1$ and $N_2$ are stationary and ergodic and
  \begin{align}
   \label{eq:condition-consistency-s}
   \esp \left[ \left(\sum_{k=1}^{N_i(1)} (\mu_i+e_{i,k}+\eta_{i,k})
     \right)^2 \right] < \infty \; ,
  \end{align}
  there exists $\varsigma>0$ such that $s_n^2 \plim \varsigma^2$.
\end{lemma}
Since the efficient shocks have zero mean and finite variance and are
independent of the point processes, a sufficient condition
for~(\ref{eq:condition-consistency-s}) to hold is
$\esp[N_i^2(1)]<\infty$ and
\begin{align*}
  \esp \left[ \left(\sum_{k=1}^{N_i(1)} \eta_{i,k} \right)^2 \right] <
  \infty \; .
\end{align*}

If $\gamma > 1/2$ and $\mu_1$ and $\mu_2$ not both zero, it follows from the
convergence~(\ref{eq:conweak-gamma>1/2}) and Lemma~\ref{lem:consistency-s} that
if the null hypothesis is true, then $t_n=O_p (n^{\gamma-1/2})$ and the
$t$-statistic diverges under the null hypothesis. Examples where this scenario
would occur include durations generated by an LMSD model with long memory and an
exponential volatility function. This scenario therefore is consistent with the
empirical properties of durations found in \citet{deo:hsieh:hurvich:2010}.

If $\gamma=1/2$ or if $\mu_1=\mu_2=0$, then $t_n \dlim \limproc(1)/\varsigma$
where $\limproc$ is the limiting process in~(\ref{eq:conweak-1/2}), i.e.~the
$t$-statistic converges to a non degenerate distribution, but the $t$-test will
not be asymptotically correctly sized, except under very specific
circumstances. First the limits $\limproc_1,\limproc_2$ must have normal
distributions. This would happen, for example, if the durations are i.i.d.~with
finite variance (as would be the case for the Poisson process), or if the
durations obey an ACD model with finite variance. Then it would also be
necessary that $\lim_{t\rightarrow \infty} t^{-1/2} \var [N_i(t)] =
\var[N_i(1)]$, which would hold for instance if $N$ is a Poisson process but
would fail if counts are autocorrelated as would typically be the case.

\subsection{A ratio statistic}
Instead of using the $t$-statistic which could be degenerate, unless the
parameter $\gamma$ is equal to 1/2, we propose to use another self-normalized
statistic which will always have a non degenerate distribution and does not need
additional assumptions. Define $\bar{r}(n) = \bar{r}_n$ and
\begin{align}
  \label{eq:ratio}
  \newstat_n = \frac{\bar{r}(n) -\mu_0^*} {|\bar{r}(n)-\bar{r}(n/2)|} \; .
\end{align}
\begin{theorem}
  \label{theo:ratio}
  If Assumptions~\ref{hypo:weakconv-counts} and~\ref{hypo:micro-negligible}
  hold, then under the null hypothesis $\mu^*=\mu_0^*$,
  \begin{align}
    \label{eq:conv-newstat}
    \newstat_n \dlim \frac{\limproc(1)}{|\limproc(1)-2\limproc(1/2)|} \; ,
  \end{align}
  where $\limproc$ denotes the limiting process in~(\ref{eq:conweak-gamma>1/2}) if
  $\gamma>1/2$ and $\mu_1$ and $\mu_2$ not both zero, or~(\ref{eq:conweak-1/2})
  if $\gamma=1/2$.  If $\mu^*>\mu_0^*$, then $\newstat_n \plim +\infty$ and if
  $\mu^*<\mu_0^*$, then $\newstat_n \plim -\infty$.
\end{theorem}

\section{Subsampling}
\label{sec:Subsampling}

Let $D_n(x)$ be the cumulative distribution of $\newstat_n$, that is, $D_n(x)=\mathbb{P}(\newstat_n \leq x)$. Let $D$ denote the corresponding limiting distribution,
\[ D = \frac{\limproc(1)}{|\limproc(1)-2\limproc(1/2)|} \]
To perform a hypothesis test on $\mu^*$ based on $T_n$, we need to approximate $D_n(x)$. Since the limiting distribution D will be unknown in practice, we will approximate $D_n(x)$ non-parametrically via subsampling.

In subsampling we split the sample into overlapping blocks of size $b$, where $b$ depends on $n$ ($b \rightarrow \infty, b/n \rightarrow 0$), given by $\{r_t, r_{t+1}, ..., r_{t+b-1}\}$, $t=1,2,...,n-b+1$, and calculate the self-normalized statistic upon each block, treating each block as if it were a full sample. Moreover, the parameter $\mu_0^*$ is replaced by its full-sample estimate $\bar{r}(n)$. This leads to $n-b+1$ subsampling statistics
\[ \newstat_{n,b,t} = \frac{\bar{r}(n,b,t) -\bar{r}(n)} {|\bar{r}(n,b,t)-\bar{r}(n,b/2,t)|},\]
\[ t=1,2,...,n-b+1; \; b \; is \; even,\]
where
\begin{align*}
\bar{r}(n)&= \frac{1}{n}\sum\limits_{k=1}^{n}r_k,\\
\bar{r}(n,b,t)&= \frac{1}{b}\sum\limits_{k=t}^{t+b-1}r_k, \\
\bar{r}(n,b/2,t) &= \frac{2}{b}\sum\limits_{k=t}^{t+b/2-1}r_k.
\end{align*}
The cumulative distribution $D_n(x)$ is then approximated by the empirical cdf of these statistics
\[ \hat{D}_{n,b}(x)=\frac{1}{n-b+1}\sum\limits_{t=1}^{n-b+1}1_{\{\newstat_{n,b,t}\leq x\}} \;. \]
Let $c_{n,b}(\alpha)$ denote the $\alpha$ quantile of the subsampling distribution, $\hat{D}_{n,b}(x)$,
\[ c_{n,b}(\alpha)=\inf\{x: \hat{D}_{n,b}(x)\geq \alpha\}. \]
Then the rejection region for the subsmapling test based on $T_n$ in the two-sided case
\begin{align*}
H_0: \mu^*=\mu_0^* \\
H_1: \mu^* \neq \mu_0^*
\end{align*}
is
\[ \{\newstat_n < c_{n,b}(\alpha/2) \;or \; \newstat_n > c_{n,b}(1-\alpha/2)\} \; .\]

To show that the resulting subsampling test based on $T_n$ has asymptotically the nominal size
($\alpha$), we use Theorem 4 of Appendix B of \citet{jach:mcelroy:politis:2012}. In order to apply
this theorem, we must establish the $\theta$-weak dependence of the time series $\{r_k\}$ (see
\citet{doukhan:louhichi:1999}; \citet{bardet:doukhan:lang:ragache:2008}). We will establish this
$\theta$-weak dependence for one of the processes we consider in Section \ref{sec:examples} (see
Proposition \ref{prop:subConsist-cox}). We say that the subsampling estimator $\hat{D}_{n,b}(x)$ is
consistent if $|\hat{D}_{n,b}(x)- D_n(x)|\plim 0$ as $n\rightarrow \infty$ for all $x$. The
following theorem can be proved based on Theorem 4, Appendix A of \citet{jach:mcelroy:politis:2012}
and its proof, as well as the discussion on page 941 of \cite{mcElroy:Jach:2012}.

\begin{theorem}
  \label{theo:subConsist}
  Assume that the time series of calendar-time returns $\{r_n,n\geq1\}$ is strictly stationary and
  $\theta$-weak dependent with rate $\theta_h = O(h^{-a})$ up to slowly varying functions. Assume
  also that the cumulative distribution function of
  \[ 
  \frac{\limproc(1)}{|\limproc(1)-2\limproc(1/2)|} 
  \] 
  in Theorem \ref{theo:ratio} is continuous. Then
  \begin{enumerate}[(i)]
  \item for $a \geq 1/2$, the subsampling estimator is consistent and the resulting test based on
    $T_n$ is asymptotically correctly sized and has power tending to 1 under the alternative
    hypothesis for any choice of the block size $b$ such that $b \rightarrow \infty$ and $b/n
    \rightarrow 0$;
  \item for $a < 1/2$, the subsampling estimator is consistent and the resulting test based on $T_n$
    is asymptotically correctly sized and has power tending to 1 under the alternative hypothesis if
    $b \rightarrow \infty$ and $b=O(n^\zeta)$ (up to slowly varying functions) for some
    $0<\zeta<2a$.
  \end{enumerate}
\end{theorem}

The theorem gives little guidance as to the choice of the block size $b$ in practice. In the
simulations of Section \ref{sec:simulation} , we will consider a variety of choices of the block
size.
% The disadvantage of this latter scenario is that $a$ is unknown to us. But in practice, block size
% will be selected by a data-driven technique as shown in \citet{jach:mcelroy:politis:2012}, that
% automatically determines on smaller blocks when serial dependence is greater.

\section{Examples}
\label{sec:examples}
In this section we give several examples of point processes which satisfy
Assumption~\ref{hypo:weakconv-counts} with~$\gamma>1/2$.

\subsection{ACD durations}
  \label{xmpl:acd}
  Assume that under the durations form an ACD(1,1) process, defined by
  \begin{align}
   \label{eq:acd}
   \tau_k & = \psi_k\epsilon_k,\qquad
   \psi_k=\omega+\alpha\tau_{k-1}+\beta\psi_{k-1},\qquad k\in\Zset,
\end{align}
where $\omega>0$ and $\alpha,\beta\geq 0$, $\{\epsilon_k\}_{k=-\infty}^\infty$
is an i.i.d.~sequence with $\epsilon_k \geq 0$ and $\esp[\epsilon_0]=1$. If
$\alpha+\beta<1$, there exists a strictly stationary solution determined by
$\tau_k = \omega \epsilon_k \sum_{j=1}^\infty \prod_{i=1}^{j-1}
(\alpha\epsilon_{k-i}+\beta)$, with finite mean
$\esp[\tau_0]=\omega/(1-\alpha-\beta)$.  The tail index $\kappa$ of a ACD(1,1)
process is the solution of the equation
\begin{align*}
  \esp[(\alpha \epsilon_0+\beta)^{\kappa}] = 1 \; .
\end{align*}
Moreover the stationary distribution satisfies $\pr(\tau_1>x) \sim c x^{-\kappa}$ for some positive
constant~$c$. See e.g.  \citet{basrak:davis:mikosch:2002}.
\begin{itemize}
\item If $1<\kappa<2$, and if $\epsilon_0$ has a positive density on $[0,\infty)$, then the finite
  dimensional distributions of $n^{-1/\kappa} \sum_{k=1}^{[nt]} (\tau_k-\esp[\tau_0])$ converges to
  a totally skewed to the right $\kappa$-stable law. Cf.
  \citet[Proposition~5]{bartkiewicz:jakubowski:mikosch:wintenberger:2011}.
  % For functional convergence in the M1 topology, see
  % \citet{basrak:krizmanic:segers:2012}.
\item A necessary and sufficient condition for $\esp[\tau_0^2]<\infty$ is
  $\esp[(\alpha\epsilon_0+\beta)^2] = \alpha^2 \esp[\epsilon_0^2]+2\alpha\beta+\beta^2<1$. Cf.
  \citet[Example~3.3]{giraitis:surgailis:2002}. Under this condition, it also holds that
  $\sum_{k=1}^\infty \cov(\tau_0,\tau_k)< \infty$.  Since the ACD process is associated, the
  summability of the covariance function implies the functional central limit theorem for the
  partial sum process. See \citet{newman:wright:1981}. Thus the
  convergence~(\ref{eq:weakconv-counts}) holds with $\gamma=1/2$ and the limit process is the
  Brownian motion.
\end{itemize}
This implies that the convergence~(\ref{eq:weakconv-counts}) holds with $\gamma=1/\kappa$ and the
limit process is a L\'evy-stable process if $\kappa<2$ and with $\gamma=1/2$ and the limit process
is the Brownian motion if $\esp[(\alpha\epsilon_0+\beta)^2] <1$.

\subsection{LMSD durations}
\label{xmpl:lmsd}
Assume that the durations form an LMSD process, defined by $\tau_k=\epsilon_k
\sigma(Y_k)$, where $\{\epsilon_{k}, k\in\Zset\}$ is an i.i.d.~sequence of
almost surely positive random variables with finite mean and $\{Y_{k},
k\in\Zset\}$ is a stationary standard Gaussian process, independent of
$\{\epsilon_{k}\}$ and $\sigma$ is a positive function. As in
\citet{deo:hsieh:hurvich:2010}, for simplicity we will assume that $\sigma(x) =
\rme^x$.

Assume that the covariance of the Gaussian process $\{Y_k\}$ is such that
\begin{align*}
  \rho_n = \cov(Y_0,Y_n) \sim c n^{2H-2} \; ,
\end{align*}
where $H \in (1/2,1)$ and $c>0$. Denote $\lambda^{-1} = \esp[\epsilon_0]
\esp[\exp(Y_0)]$. Then we have the following possibilities.
\begin{itemize}
\item If $\esp[\epsilon_k^2]<\infty$, then
  \begin{align*}
    n^{-H} \sum_{k=1}^{[n\cdot]} (\tau_k-\lambda^{-1}) \func \varsigma B_{H} \; ,
  \end{align*}
  where $\varsigma$ is a nonzero constant and $B_{H}$ is the standard fractional
  Brownian motion.
\item If $\pr(\epsilon_1>x) \sim c x^{-\alpha}$ as $x\to\infty$ with
  $\alpha\in(1,2)$, then $\esp[\epsilon_k^2]=\infty$ and the following dichotomy
  is proved in \citet{kulik:soulier:2012}.
  \begin{itemize}
  \item If $H > 1 - 1/\alpha$, then
  \begin{align*}
    n^{-H} \sum_{k=1}^{[n\cdot]} (\tau_k-\lambda^{-1}) \func \varsigma B_{H} \;    .
  \end{align*}
\item If $H < 1 - 1/\alpha$, then
  \begin{align*}
    n^{-1/\alpha} \sum_{k=1}^{[n\cdot]} (\tau_k-\lambda^{-1}) \func L_\alpha \; ,
  \end{align*}
  where $L_\alpha$ is a totally skewed to the right $\alpha$-stable L\'evy process.
\end{itemize}
\end{itemize}
We thus see that~(\ref{eq:weakconv-counts}) may hold with $\gamma=H$ in the first case and $\gamma =
1/\alpha$ with a stable non Gaussian limit in the latter case.

\subsection{Superposition of independent point processes}

This is an example showing that the limiting distribution of a centered and normalized counting
process $N$ may be stable (hence heavy tailed) even though its durations $\tau_k$ are light tailed.
Let $M_1$ be a Poisson process with intensity $\lambda$ and durations $\{\tau_k^{(1)}\}$ and $M_2$
be a renewal process, independent of $M_1$, with i.i.d. durations~$\{\tau_k^{(2)}\}$ in the normal
domain of attraction of a stable law, i.e.
\begin{align*}
  n^{-1/\alpha} \sum\limits_{k=1}^{[nt]} (\tau_k^{(2)}-\lambda_2^{-1})
  \fidi \mathbf{S}(t) \; ,
\end{align*}
where $\alpha\in(1,2)$, $\lambda_2^{-1} = \esp[\tau_1^{(2)}]$ and $\mathbf S$ is a totally skewed to
the right $\alpha$-stable L\'evy process. By the CLT equivalence, this implies that
\begin{align*}
  n^{-1/\alpha} (M_2(nt) -\lambda_2 nt) \fidi  \mathbf{S}(t) \; .
\end{align*}
Since the durations of the Poisson process are i.i.d. with finite variance, $n^{-1/2}(M_1(nt) -
\lambda_1nt)$ converges weakly to a Brownian motion. Thus defining the superposition $N$ of these
two point processes by
\begin{align*}
  N = M_1+M_2 \; ,
\end{align*}
we obtain that
\begin{align*}
  n^{-1/\alpha} (N(nt) -\lambda nt) \fidi -\lambda_2^{1+1/\alpha} \mathbf{S}(t) \; ,
\end{align*}
with $\lambda=\lambda_1+\lambda_2$.

Let $\{\tau_k\}$ be the duration sequence of the point process $N$.  We now show that the durations
are light tailed under the Palm measure $P^0$ under which they form a stationary sequence.  It
follows from the formula (1.4.5) and example 1.4.1 in \citet{baccelli:bremaud:2003} that
\begin{align*}
  P^0(\tau_1>x) & = \lambda^{-1} \rme^{-\lambda_1 x} \left\{ \lambda_1 \bar F_2(x) + \lambda_2 \bar
    H_2(x) \right\} \; ,
\end{align*}
where $F_2$ is the distribution function of $\tau_1^{(2)}$ and $H_2$ is the corresponding delay
distribution defined by $H_2(x) = \lambda \int_0^x \bar F_2(t) \, \rmd t$. This yields that
$P^0(\tau_1>x) \leq \rme^{-\lambda_1 x}$, i.e. the durations of the superposition process are light
tailed.

\subsection{Cox processes}
Consider now a Cox process $N$ driven by a stationary random measure $\xi$, which means that
conditionally on $\xi$, $N$ is a Poisson point process with mean measure $\xi$. Then, denoting
$\xi(t)=\xi([0,t])$, we have $\esp[N(t)] = \esp[\xi(t)] = \lambda t$ with $\lambda = \esp[\xi(1)]$
and
\begin{align*}
  \var(N(t)) = \esp[\xi(t)] + \var(\xi(t)) = \lambda t + \var(\xi(t)) \; .
\end{align*}
If the stationary random measure $\xi$ has long memory with Hurst index $H>1/2$ then so has $N$. We
give examples.
\begin{enumerate}[(i)]
\item Consider the random measure $\xi$ with density $\rme^{\sigma Z_H(s)}$ with respect to
  Lebesgue's measure on $\Rset$, where $\sigma>0$ and $Z_H$ is a standard fractional Gaussian noise
  with Hurst index $H\in(1/2,1)$, i.e. a Gaussian stationary process with covariance function
  $\rho_H(t) = \frac12\{|t-1|^{2H}-2|t|^{2H}+|t+1|^{2H}\}$. This means that for any Borel set $A$,
  $\xi(A) = \int_A \rme^{Z_H(s)} \, \rmd s$, and in particular $\xi(t) = \int_0^t \rme^{Z_H(s)}\,
  \rmd s$. Then, by stationarity,
  \begin{align*}
    \lambda = \esp \left [ \int_0^1\rme^{\sigma Z_H(s)} \, \rmd s \right] = \esp
    \left[ \rme^{\sigma Z_H(0)} \right] = \rme^{-\sigma^2/2} \; .
  \end{align*}

  Since the function $x\to\rme^x$ has Hermite rank~1, it holds that
  $\cov(\rme^{Z_H(0)},\rme^{Z_H(t)}) \sim c\rho_H(t)$ and $\var(\xi(t)) \sim C^2 t^{2H}$ where $c$
  and $C$ are positive constants, and $n^{-H} \{\int_0^{nt} \rme ^{Z_H(s)} \, \rmd s - \lambda nt \}
  \fidi C B_H(t)$. See \citet{arcones:1994} and~\citet{dobrushin:major:1979} for details about
  Hermite ranks, covariance inequality and convergence in distribution of integrals of functions of
  Gaussian processes.  Applying Lemma~\ref{lem:conv-cox}, we obtain the convergence
  \begin{align*}
    n^{-H} \{N(nt)  - \lambda nt \} \fidi C B_H(t) \; .
  \end{align*}

\item Consider the case where the stochastic intensity can be expressed as $\xi(t) = \int_0^t W(s)
  \, \rmd s$, where $W(t)$ is an alternating renewal or ON-OFF process
  (cf.\citet{heath:resnick:samorodnitsky:1998}, \citet{daley:2010}), defined by $W(t)=1$ if
  $T_{2n}\leq t <T_{2n+1}$ and $T_{2n+1} \leq t < T_{2n+2}$, where $T_{2n}=\sum_{i=0}^n X_i+Y_i$ and
  $T_{2n}=\sum_{i=0}^n X_i+Y_i+X_{n+1}$, $\{X_i,i\geq1\}$ and $\{Y_i,i\geq1\}$ are two independent
  i.i.d. sequences of positive random variables with finite mean (the ON and OFF periods) and $X_0$
  and $Y_0$ are mutually independent and independent of the other random variables, whose
  distributions are the delay distributions which make the renewal processes stationary. If the ON
  distribution is regularly varying at infinity with tail index $\alpha\in(1,2)$
  i.e. $\pr(X_1>x)\sim c x^{-\alpha}$ as $x\to\infty$ (which implies infinite variance) and if
  $\esp[Y_i^{\alpha+\epsilon}]<\infty$ for some $\epsilon>0$, then $\var(\xi(t))\sim Ct^{2H}$ with
  $H=(3-\alpha)/2$. Moreover, $n^{-1/\alpha} \{\xi(nt) - \lambda nt\} \fidi \Lambda(t)$, where
  $\Lambda$ is a totally skewed to the right $\alpha$-stable process. Applying
  Lemma~\ref{lem:conv-cox}, we obtain the convergence
  \begin{align*}
    n^{-1/\alpha} \{N(nt)  - \lambda nt \} \fidi \Lambda(t) \; .
  \end{align*}

\end{enumerate}

\begin{proposition}
  \label{prop:subConsist-cox}
  Consider the returns $\{r_n\}_{n \in \mathbb{N}^*}$ driven by a Cox process of case (i), then the
  conclusion of Theorem \ref{theo:subConsist} holds with $a=2-2H$.
\end{proposition}

\begin{proof}[Proof of Proposition~\ref{prop:subConsist-cox}]
  First note that the causal Gaussian process $Z_H(s)$ is $\theta$-weak dependent with rate
  $\theta_h = L(h)h^{2H-2}$ by \citet{bardet:doukhan:lang:ragache:2008}, where $L(h)$ is a slowly
  varying function. Then by the Proposition 1 of \citet{jach:mcelroy:politis:2012}, the stationary
  process $w(s) = e^{\sigma Z_H(s)}$ is $\theta$-weak dependent with the same rate. Using
  Lemma~\ref{lem:WeakDependenceForCox}, we have that the returns $\{r_n\}_{n \in \mathbb{N}^*}$ are
  $\theta$-weak dependent with the same rate $\theta_h = O(h^{2H-2})$. The cumulative distribution
  function of
  \[
  \frac{\Lambda(1)}{|\Lambda(1)-2\Lambda(1/2)|} 
  \]
  is continuous, where $\Lambda(1)$ and $\Lambda(1/2)$ have totally skewed to the right
  $\alpha$-stable distributions. Finally, by applying Theorem~\ref{theo:subConsist} with $a=2-2H$,
  the result is proved.
\end{proof}

\section{Simulation}
\label{sec:simulation}

We study the performance of the ordinary $t$-test and the subsampling-based $T$ test of
\begin{align}
  \label{eq:hypo}
  &H_0: \mu^* = 0\\
  \label{eq:hypoalt}
  &H_1: \mu^* > 0
\end{align}

We consider exponential, ACD(1,1), and LMSD durations. In each case, we simulate two mutually
independent durations process $\{\tau_{k}^{(i)}\}$ for $i=1,2$, and then obtain the corresponding
counting processes $\{N_i(t)\}$. Next, we generate mutually independent disturbance series $\{\tilde
e_{k}^{(1)}\}$ and $\{\tilde e_{k}^{(2)}\}$ that are i.i.d. Gaussian with means $\mu_1>0$,
$\mu_2<0$, $\mu_1+\mu_2>0$, and with finite variances $\sigma_i^2$. For simplicity, we assume all
parameters other than $\mu_i$ for the two processes to be the same. For the exponential and
EACD(1,1) models, we only include the efficient shocks in simulations, and for the LMSD, we also
include microstructure shocks, $\{\eta_k^{(i)}\}$, that induce a leverage effect. We then construct
the log-price series $\{y(k)\}_{k=1}^n$ and return series $\{r_k\}_{k=1}^n$ from (\ref{eq:trend})
and $r_k=y(k)-y(k-1)$ with $y(0)=0$. We calculate $t_n$ and $T_n$ based on return series and apply
subsampling to obtain the empirical quantiles of $T_n$ and thereby carry out the hypothesis test.

We do a preliminary study on the sensitivity of our results to the choice of the block size
$b$. These results show insensitivity except for the largest values of $b$.

The simulation results show that in all three models, the $T$ test is generally correctly sized. But
in the exponential case where the usual $t$ test is also correctly sized, the power of $T$ is lower
than that of the $t$ test. To improve the power, we modify the test statistic defined in
(\ref{eq:ratio}) as
\begin{align}
  \label{eq:newratio}
  \newstat_{n,2} = \frac{\bar{r}(n,n,1) -\mu_0^*} {|\bar{r}(n,n/2,1)-\bar{r}(n,n/4,1)| +
    |\bar{r}(n,n/2,n/2+1)-\bar{r}(n,n/4,n/2+1)|} \; ,
\end{align}
where $\bar{r}(n,b,t)= \frac{1}{b}\sum\limits_{k=t}^{t+b-1}r_k$. In the remainder of this section,
we will use $T_1$ to denote (\ref{eq:ratio}), $T_2$ to denote (\ref{eq:newratio}), and ``$T$ test"
to denote a test based on either $T_1$ or $T_2$.

\subsection{Exponential durations}
\label{subsec:expsimu}

We generate $\{\tau_{k}^{(1)}\}$, $\{\tau_{k}^{(2)}\} $ from exponential distributions
$\tau_{k}^{(i)}\iid Exp(\lambda_i)$, which result in $N_1(t)$ and $N_2(t)$ being Poisson
processes. Since the returns are i.i.d, the $t$-test is correctly sized. Our goal is to study the
size of the $T$ test and to compare the power of the $T$ to that of the $t$ test.

\subsubsection{Parameter calibration}
\label{subsec:expParCali}

We consider a total of 30 configurations of parameters and sample size, 15 for evaluating size
(Table \ref{tab:poissize}), and 15 for evaluating power (Table \ref{tab:poispower}). The model has
six parameters ($\lambda_1$, $\lambda_2$, $\sigma_1$, $\sigma_2$, $\mu_1$, $\mu_2$). Note that we
assume that ($\lambda_1 = \lambda_2=\lambda$, $\sigma_1 = \sigma_2=\sigma$). For the first
configuration corresponding to power (Table \ref{tab:poispower}), we calibrate the model to match
the mean and standard deviation of the returns to those of the Fama/French factor, Rm-Rf (see
Section \ref{sec:data}), and the intensity observed in the Boeing series used in
\citet{deo:hsieh:hurvich:2010}.

We consider one unit of time to represent 5 minutes, whereas \citet{deo:hsieh:hurvich:2010} used 1
minute and Rm-Rf is a daily series. Our calibrated model (which accounts for the discrepancies in
the time units) is given in the first row of Table \ref{tab:poispower}. We compute the parameter
values using the fact that for the Poisson process, there is an explicit functional relationship
between the model parameters and the mean and variance of the calendar-time returns.

\begin{figure}[H]
\begin{center}
  \includegraphics[scale=.8, trim = 30mm 145mm 30mm 40mm]{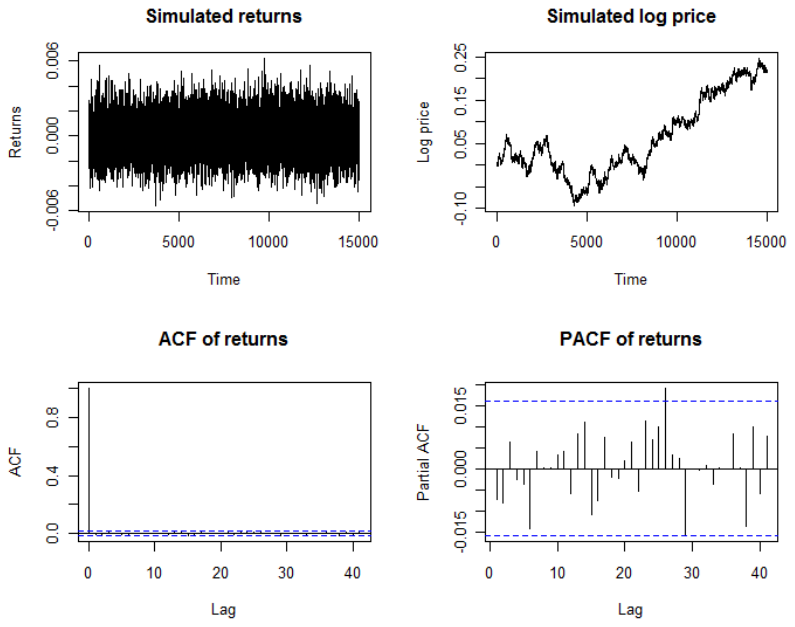}
  \caption{\footnotesize One simulated realization of the calibrated Poisson model (first row in
    Table \ref{tab:poispower}). Time series plot of returns and log price, as well as ACF and PACF,
    $n=15000$. The model is calibrated using Boeing's intensity and Fama/French factor Rm-Rf, such
    that $\lambda=2.71$, $\sigma=0.0003$, and the mean and standard deviation for the 5-minute
    returns are $\mu^*= 3.846\times 10^{-6}$, $\tilde{\sigma} = 0.0012$. }\label{fig:simuPois}
\end{center}
\end{figure}

Figure~\ref{fig:simuPois}, for one simulated realization from this model with $n=15000$, shows the
time series plots of returns and log price, as well as the ACF and PACF of the returns. Neither ACF
nor PACF shows statistically significant lags.

\begin{figure}[H]
\begin{center}
  \includegraphics[scale=.73, trim = 45mm 130mm 30mm 40mm]{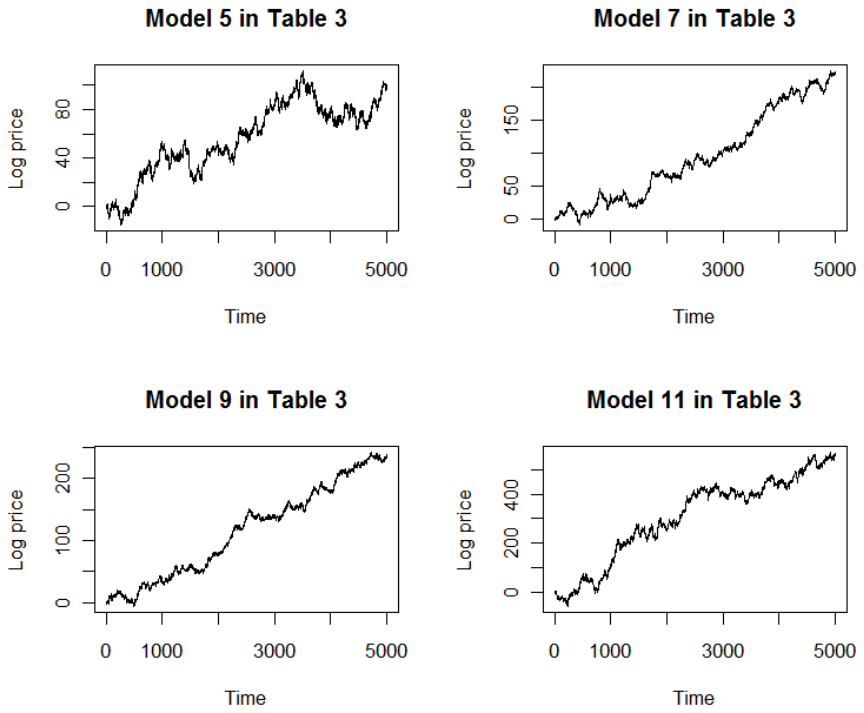}
  \caption{\footnotesize Simulated log price for models in rows 5 ($\mu^*=3.846 \times 10^{-3}$), 7,
    9, and 11 (all with $\mu^*=3.846 \times 10^{-1}$) of Table \ref{tab:poispower},
    $n$=5000.}\label{fig:simuPois2}
\end{center}
\end{figure}

For the remaining configurations in Table \ref{tab:poispower}, we increase $\mu^*$ by varying
$\lambda$ and $\mu_i$. We also consider the two sample sizes, $n=5000$ (corresponding to roughly one
quarter of 5-minute observations) and $n=10000$. The power should increase as $\mu^*$
increases. Figure \ref{fig:simuPois2} shows the time series plots of the log price for the models in
rows 5 ($\mu^*=3.846 \times 10^{-3}$), 7, 9, and 11 (all with $\mu^*=3.846 \times 10^{-1}$) in Table
\ref{tab:poispower}.

For the first configuration corresponding to size (Table \ref{tab:poissize}), we only match the
standard deviation of the returns to that of Rm-Rf and the intensity of Boeing series, while setting
the mean to zero. For the remaining configurations in Table \ref{tab:poissize}, we vary the
$\lambda$ and $\mu_i$ similarly as in Table \ref{tab:poispower} while keeping $\mu^*=0$.

For each configuration in Tables \ref{tab:poissize} and \ref{tab:poispower}, we generate 1000
realizations. For each realization, we calculate the $t$ and $T$ statistics, and use subsampling to
obtain the rejection region of $T$ for the hypothesis test (\ref{eq:hypo}) (\ref{eq:hypoalt}) at the
0.05 significance level,
\[ 
\{T > c_{n,b}(0.95)\} 
\]
where $c_{n,b}(p)$ is the $p^{th}$ quantile of the subsampling distribution as defined in Section
\ref{sec:Subsampling}.

\subsubsection{Block size}
\label{subsec:expBlocksize}

For the two calibrated models, we evaluate the size and power of $T_1$ with different block sizes
$b$, sample size $n=15000$ (corresponding to roughly three quarters of observations), and 1000
realizations. We also evaluate the size and power. Results are given in Table
\ref{tab:poisTsubsize}. When $b<1600$, the size of $T_1$ is not significantly different from 0.05
using the binomial test at the 0.05 significance level. When $b$ becomes larger, the size is
significantly larger than 0.05. The size and power of the $T_1$ test are insensitive to the block
size except at the largest values of $b$. In view of this insensitivity, and in view of the fact that in Theorem \ref{theo:subConsist} $b$ may approach $\infty$ arbitrarily slowly, we will fix the block size
$b=160$ in the following simulations. (Other values of $b$ were tried there and gave similar results.)

\begin{table}[H]
	\begin{center}
		\footnotesize
		\begin{tabular}{*{9}{c}|c}
			\hline
			\multicolumn{9}{c|}{\textbf{$T_1$ test}} & \textbf{$t$ test} \\\hline
			\textbf{Block size ($b$)}& 160 & 320& 480& 640& 800& 1120& 1600 & 3200 &\\
			\textbf{Size}&0.053&0.052&0.053&0.053&0.053&0.054&0.065&0.077&0.053\\
			\textbf{Power}	&0.071&0.07&0.071&0.07&0.078&0.082&0.090&0.110&0.089\\\hline
		\end{tabular}
		\caption{\footnotesize Size and power of $T_1$ with different block sizes,
                  $n=15000$, exponential durations, using the calibrated models (from first row of
                  Tables \ref{tab:poissize} and \ref{tab:poispower}).} \label{tab:poisTsubsize}
	\end{center}
\end{table}

\subsubsection{Size/Power of the $t$ and $T$ tests for exponential durations}
\label{subsec:expSizePower}

Table \ref{tab:poissize} shows that both $t$ and $T$ are correctly sized. Indeed, none of the sizes
reported in Table \ref{tab:poissize} are significantly different from 0.05 according to the binomial
test at the 0.05 significance level. Table \ref{tab:poispower} shows that in general the power of
$T$ is lower than that of $t$, and power of all tests increases as $\mu^*$ increases. Furthermore,
as $n$ is increased holding all model parameters fixed, the power of all tests increases. Finally,
the power of the $T$ tests is lower than that of the $t$ test, but $T_2$ has higher power than
$T_1$.

%

%\begin{landscape}
\begin{table}[H]
\begin{center}
\small
\begin{tabular}{ccccc|ccc}
\hline
 \multicolumn{5}{c|}{\textbf{Parameter values}} &\multicolumn{3}{c}{\textbf{Size}}  \\
$n$ & $\lambda_1=\lambda_2$ & $\mu_1$ & $\mu_2$ & $\mu^*$  & $T_1$ & $T_2$ & $t$ \\ \hline
5000& 2.71& $4.266 \times 10^{-4}$ & $-4.266 \times 10^{-4}$ & 0 & 0.055  & 0.063 & 0.045 \\
5000& 2.71& $4.173 \times 10^{-4}$ & $-4.173 \times 10^{-4}$ & 0 & 0.057  & 0.050 & 0.048\\
5000& 2.71& $4.359 \times 10^{-4}$ & $-4.359 \times 10^{-4}$ & 0 & 0.048  & 0.056 & 0.048 \\
5000& 27.1& $4.266 \times 10^{-4}$ & $-4.266 \times 10^{-4}$ & 0 & 0.055  & 0.054 & 0.052 \\
5000& 2.71& 0.4266 & -0.4266 & 0 & 0.043 & 0.054&0.049\\
5000& 2.71& 0.4173 & -0.4173 & 0 & 0.047 & 0.055&0.049\\
5000& 2.71& 0.3273 & -0.3273 & 0 & 0.058 & 0.056&0.053\\
5000& 2.71& 0.4359 & -0.4359 & 0 & 0.061 & 0.057&0.063\\
5000& 2.71& 0.5259 & -0.5259 & 0 & 0.060 & 0.058&0.053\\
5000& 27.1& 0.4266 & -0.4266 & 0 & 0.058 & 0.057&0.046 \\
5000& 271 & 0.4266 & -0.4266 & 0 & 0.053 & 0.051&0.058 \\
10000&2.71& 0.4173 & -0.4173 & 0 & 0.050 & 0.064&0.056\\
10000&2.71& 0.3273 & -0.3273 & 0 & 0.052 & 0.051&0.050\\
10000&2.71& 0.4359 & -0.4359 & 0 & 0.058 & 0.056&0.057\\
10000&2.71& 0.5259 & -0.5259 & 0 & 0.060 & 0.058&0.056\\\hline
\end{tabular}
\caption{\small Size of the $t$ and $T$ tests for i.i.d. exponential durations.} \label{tab:poissize}
\end{center}
\end{table}
%\end{landscape}

%\begin{landscape}
\begin{table}[H]
\begin{center}
\small
\begin{tabular}{ccccc|ccc}
\hline
 \multicolumn{5}{c|}{\textbf{Parameter values}} &\multicolumn{3}{c}{\textbf{Power}}  \\
$n$ & $\lambda_1=\lambda_2$ & $\mu_1$ & $\mu_2$ & $\mu^*$  & $T_1$ & $T_2$ & $t$ \\ \hline
5000& 2.71& $4.273 \times 10^{-4}$ & $-4.259 \times 10^{-4}$ & $3.846 \times 10^{-6}$  &0.066 &0.077 & 0.068\\
5000& 2.71& $4.273 \times 10^{-4}$ & $-4.173 \times 10^{-4}$ & $3.846 \times 10^{-5}$  &0.237 &0.363 &0.575\\
5000& 2.71& $4.359 \times 10^{-4}$ & $-4.259 \times 10^{-4}$ & $3.846 \times 10^{-5}$   &0.242 &0.242&0.582 \\
5000& 27.1& $4.273 \times 10^{-4}$ & $-4.259 \times 10^{-4}$ & $3.846 \times 10^{-5}$ &0.103&0.105&0.147 \\
5000& 2.71& 0.4273 & -0.4259 & $3.846 \times 10^{-3}$  &0.068 &0.089&0.090\\
5000& 2.71& 0.4273 & -0.4173 & $3.846 \times 10^{-2}$  &0.361 &0.615&0.878\\
5000& 2.71& 0.4273 & -0.3273 & $3.846 \times 10^{-1}$  &1 &1 &1\\
5000& 2.71& 0.4359 & -0.4259 & $3.846 \times 10^{-2}$  &0.359 &0.593&0.853\\
5000& 2.71& 0.5259 & -0.4259 & $3.846 \times 10^{-1}$ &1 &1   &1\\
5000& 27.1& 0.4273 & -0.4259 & $3.846 \times 10^{-2}$  &0.151&0.204&0.222 \\
5000& 271 & 0.4273 & -0.4259 & $3.846 \times 10^{-1}$  &0.357 &0.596&0.880 \\
10000&2.71& 0.4273 & -0.4173 & $3.846 \times 10^{-2}$  &0.492 &0.788&0.993\\
10000&2.71& 0.4273 & -0.3273 & $3.846 \times 10^{-1}$  &1 &1 &1\\
10000&2.71& 0.4359 & -0.4259 & $3.846 \times 10^{-2}$  &0.461 &0.790 &0.991\\
10000&2.71& 0.5259 & -0.4259 & $3.846 \times 10^{-1}$ &1  &1  &1\\\hline
\end{tabular}
\caption{\small Power of the $t$ and $T$ tests for i.i.d. exponential durations.} \label{tab:poispower}
\end{center}
\end{table}
%\end{landscape}

\subsection{ACD durations}
\label{subsec:acdsimu}

We generate $\{\tau_{k}^{(i)}\}$ from the exponential ACD(1,1) model (EACD(1,1)) that satisfies
(\ref{eq:acd}) with $\{\epsilon_{k}^{(i)}\} \iid exp(1)$, $i=1,2$. We generate models with finite
and infinite variance for the durations. For the finite variance model, we try two block sizes
$b=160$ and 320. For the EACD(1,1), $t$ is oversized and therefore we will not examine its power. We
will also study the size and power of the $T$ tests.

\subsubsection{Parameter calibration}
\label{subsec:acdParCali}

Similarly as in Section \ref{subsec:expsimu}, for both finite and infinite variance EACD models, we consider a total of 30 configurations of parameters and sample size, 15 for evaluating size (Table \ref{tab:eacdsizeFin} for the finite variance model, Table \ref{tab:eacdsizeInf} for the infinite variance model), and 15 for evaluating power (Table \ref{tab:eacdpowerFin} for the finite variance model, Table \ref{tab:eacdpowerInf} for the infinite variance model). The model has six parameters ($\alpha, \beta, \omega, \sigma, \mu_1, \mu_2$). Note that we assume that $\alpha, \beta, \omega, \sigma$ are the same for both processes $\{\tau_{k}^{(1)}\}$, $\{\tau_{k}^{(2)}\}$. \\

\noindent \textit{EACD(1,1) with finite variance}

The necessary and sufficient condition for finite variance of the EACD(1,1) model is
\begin{align}
\label{eq:finCond}
2\alpha^2+\beta^2+2\alpha\beta = \alpha^2 + (\alpha+\beta)^2<1 \;,
\end{align}

For the first configuration corresponding to power (Table \ref{tab:eacdpowerFin}), we still calibrate the model to match the mean and standard deviation of the returns to those of the Fama/French factor, Rm-Rf, and the ($\alpha, \beta, \omega$) observed in the Boeing series used in \citet{deo:hsieh:hurvich:2010}, which satisfy (\ref{eq:finCond}). We use the same $\sigma$ as for the exponential durations and set $\mu_1$ and $\mu_2$ to match the mean and standard deviation of the Rm-Rf. This calibrated model also accounts for the discrepancies in the time units, and is given in the first row of Table \ref{tab:eacdpowerFin}. For the EACD(1,1) model, there is no explicit functional relationship between the model parameters and the mean and variance of the calendar-time returns, since these now also depend on the variance of counts (the number of events occurring in a given 5-minute time interval). We estimate the variance of counts using simulation and set $\mu_1$ and $\mu_2$ accordingly.

\begin{figure}[H]
\begin{center}
  \includegraphics[scale=.8, trim = 30mm 145mm 30mm 20mm]{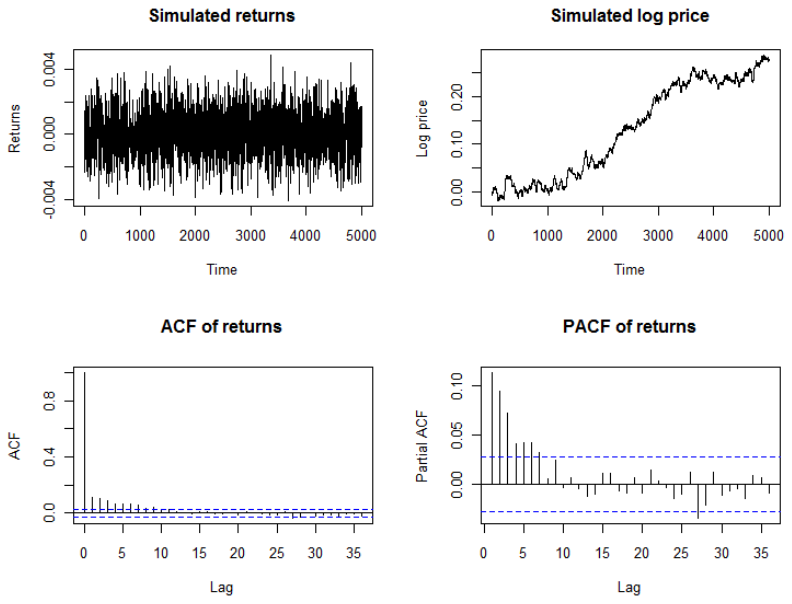}
  \caption{\footnotesize One simulated realization of the calibrated EACD(1,1) model with finite
    variance (first row in Table \ref{tab:eacdpowerFin}). Time series plot of returns and log price,
    as well as ACF and PACF, $n=5000$. The model is calibrated using Boeing's and Fama/French factor
    Rm-Rf, such that $\alpha=0.023161$ $\beta=0.970158$, $\omega=0.00247$, $\sigma=0.0003$, and the
    mean and standard deviation for the 5-minute returns are $\mu^*= 3.846\times 10^{-6}$,
    $\tilde{\sigma} = 0.0012$. }\label{fig:simuACDfin}
\end{center}
\end{figure}

Figure~\ref{fig:simuACDfin}, for one simulated realization from this model with $n=5000$, shows the time series plots of returns and log price, as well as the ACF and PACF of the returns. Both ACF and PACF show statistically significant autocorrelations at several lags.	

\begin{figure}[H]
  \begin{center}
    \includegraphics[scale=.72, trim = 45mm 130mm 30mm 40mm]{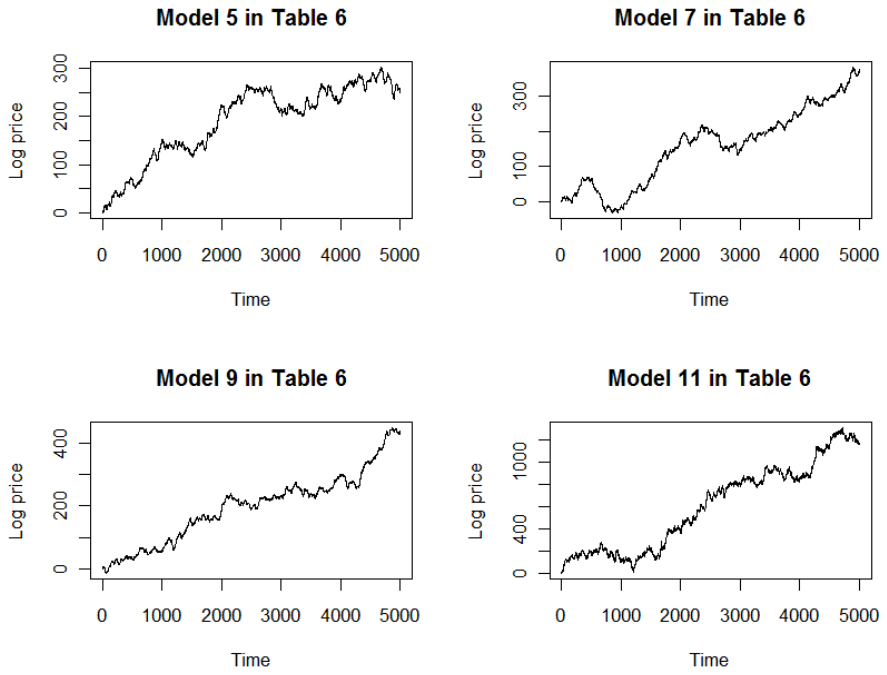}
    \caption{ \footnotesize Simulated log price for models in rows 5 ($\mu^*=3.846 \times 10^{-3}$),
      7, 9, and 11 (all with $\mu^*=3.846 \times 10^{-1}$) of Table \ref{tab:eacdpowerFin},
      $n$=5000. }\label{fig:simuACDfin2}
  \end{center}
\end{figure}

For the remaining configurations in Table \ref{tab:eacdpowerFin}, we increase $\mu^*$ by varying $\omega$ (which changes $\lambda$) and $\mu_i$. We also consider the two sample sizes, $n=5000$ and $n=10000$. The power should generally increase as $\mu^*$ increases. Figure \ref{fig:simuACDfin2} shows the time series plots of the log price for the models of rows 5 ($\mu^*=3.846 \times 10^{-3}$), 7, 9, and 11 (all with $\mu^*=3.846 \times 10^{-1}$) of Table \ref{tab:eacdpowerFin}.

For the first configuration corresponding to size (Table \ref{tab:eacdsizeFin}), similarly as in Section \ref{subsec:expParCali}, we set the mean to zero. For the remaining configurations in these two tables, we vary the $\omega$ and $\mu_i$ similarly as in Table \ref{tab:eacdpowerFin} while keeping $\mu^*=0$.

For each configuration of size and power, we generate 1000 realizations. The simulation procedure is the same as for the exponential durations. \\

\noindent \textit{EACD(1,1) with infinite variance}

The necessary and sufficient condition for the EACD(1,1) model to have infinite variance is
\[ \alpha, \beta \geq 0, \;\; \alpha+\beta<1, \;\;  \alpha^2 + (\alpha+\beta)^2\geq 1 \;.\]
We fix $\alpha+\beta=0.99$ which is the same as for all of the finite variance models. Table \ref{tab:eacdinfpara} shows the ($\alpha$, $\beta$) values that satisfy the above constraint while attaining the values of $\alpha^2 + (\alpha+\beta)^2$ given in the top row. Since for the ACD model, estimates of $\beta$ are typically close to 1, we choose the first pair $\alpha=0.1729$, $\beta=0.817$. We adjust $\omega$ to obtain the same values of $\lambda$ as used in the finite variance EACD(1,1) models. The value of $\sigma$ is the same as before. Then following the same procedure as described in the finite variance model, we obtain the size (Table \ref{tab:eacdsizeInf}) and power (Table \ref{tab:eacdpowerInf}).

\begin{figure}[H]
  \begin{center}
    \includegraphics[scale=.8, trim = 30mm 150mm 30mm 38mm]{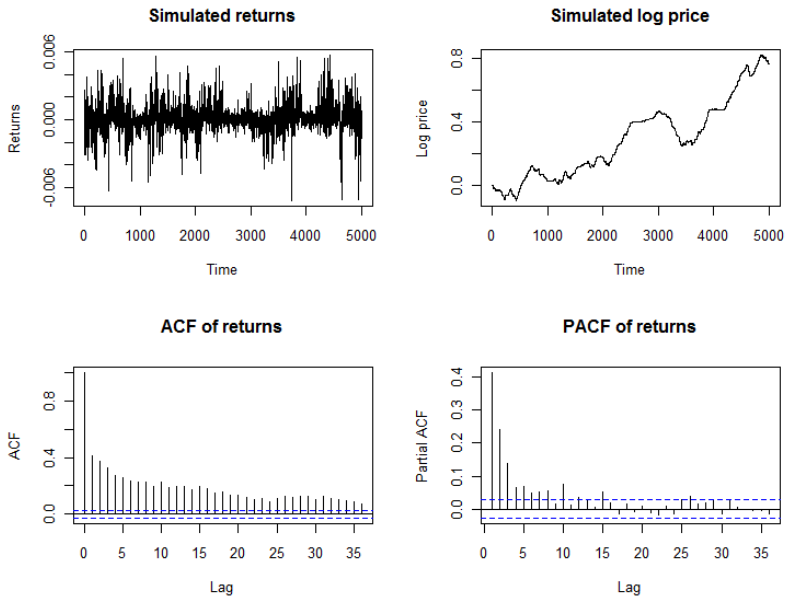}
    \caption{\footnotesize One simulated realization of the calibrated EACD(1,1) model with infinite
      variance (first row in Table \ref{tab:eacdpowerInf}). Time series plot of returns and log
      price, as well as ACF and PACF, $n=5000$. The model is calibrated using Boeing's and
      Fama/French factor Rm-Rf, such that $\alpha=0.1729$, $\beta=0.8171$, $\omega=0.00369$,
      $\sigma=0.0003$, and the mean and standard deviation for the 5-minute returns are $\mu^*=
      3.846\times 10^{-6}$, $\tilde{\sigma} = 0.0012$. }\label{fig:simuACDinf}
  \end{center}
\end{figure}

\begin{figure}[H]
  \begin{center}
    \includegraphics[scale=.72, trim = 45mm 130mm 30mm 40mm]{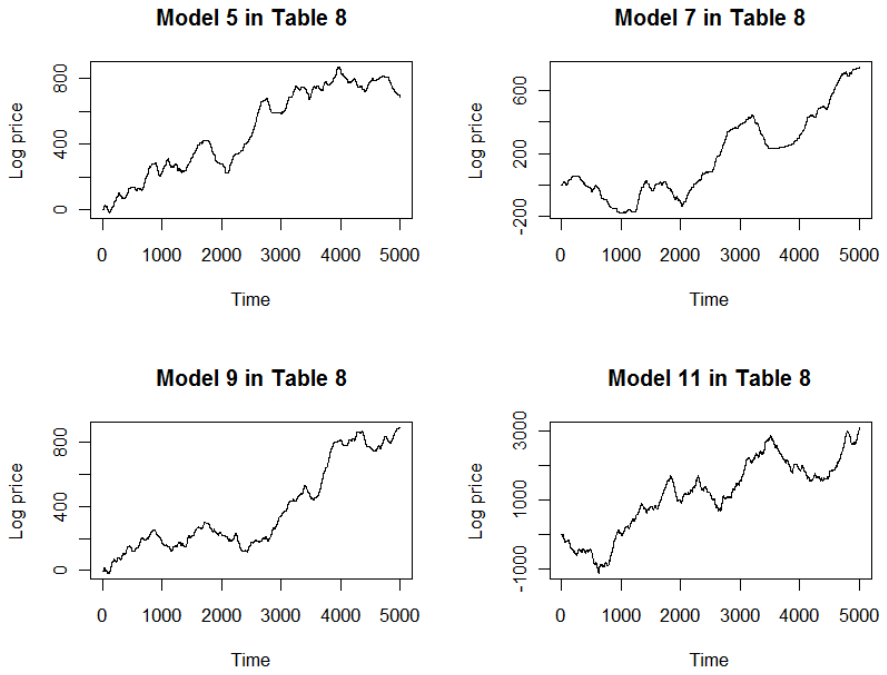}
    \caption{\footnotesize Simulated log price for models in rows 5 ($\mu^*=3.846 \times 10^{-3}$),
      7, 9, and 11 (all with $\mu^*=3.846 \times 10^{-1}$) of Table \ref{tab:eacdpowerInf},
      $n$=5000.}\label{fig:simuACDinf2}
  \end{center}
\end{figure}

Figure~\ref{fig:simuACDinf}, for one simulated realization from this model with $n=5000$, shows the time series plots of returns and log price, as well as the ACF and PACF of the returns. Figure \ref{fig:simuACDinf2} shows the time series plots of the log price for the models of rows 5 ($\mu^*=3.846 \times 10^{-3}$), 7, 9, and 11 (all with $\mu^*=3.846 \times 10^{-1}$) of Table \ref{tab:eacdpowerInf}.

We obtained similar results not shown here for other configurations of ($\alpha$, $\beta$, $\omega$)  both with and without the constraint $\alpha+\beta=0.99$.

\begin{table}[H]
	\begin{center}
		\small
		\begin{tabular}{c|cccccc}
			\hline
			$\alpha^2 + (\alpha+\beta)^2$ & 1.01 & 1.02 & 1.04 & 1.06 & 1.08 & 1.10\\\hline
			$\alpha$ & 0.1729 & 0.1997 & 0.2447 & 0.28267 & 0.3161 & 0.3463\\
			$\beta$& 0.8171&0.7903&0.7453& 0.7073 & 0.6739 & 0.6437\\\hline
		\end{tabular}
		\caption{\small ($\alpha$, $\beta$) values such that $\alpha+\beta=0.99$ and $\alpha^2 + (\alpha+\beta)^2$ is given in the top row.} \label{tab:eacdinfpara}
	\end{center}
\end{table}

\subsubsection{Size/Power of the $t$ and $T$ tests for EACD(1,1) durations}
\label{subsec:acdSizePower}

The size for both the finite variance (Table \ref{tab:eacdsizeFin}) and infinite variance models (Table \ref{tab:eacdsizeInf}) show that for EACD(1,1) durations, the $t$ test is over-sized and not reliable, especially in the infinite variance cases. Thus, we will not examine the power of the $t$ test here.

For the finite variance cases, the $T$ test is correctly sized. Table \ref{tab:eacdpowerFin} indicates that the power of $T$ has the same properties as in the Poisson model, i.e., the power increases as $\mu^*$ or $n$ increases; $T_2$ generally has higher power than $T_1$.

For the infinite variance cases, Table \ref{tab:eacdsizeInf} for size shows that both $T_1$ and $T_2$ are under-sized in some cases and are approximately correctly sized in the other cases. Table \ref{tab:eacdpowerInf} shows that the power of both $T_1$ and $T_2$ increases as $\mu^*$ or $n$ increases, but perhaps not as fast as for the finite variance models. $T_2$ has higher power than $T_1$, particularly when $\mu^*$ is large.

The results of Tables \ref{tab:eacdsizeFin} and \ref{tab:eacdpowerFin} show no apparent sensitivity to the block size.

%\begin{landscape}
\begin{table}[H]
\begin{center}
\footnotesize
\begin{tabular}{cccccc|cc|cc|c}
\hline
\multicolumn{6}{c|}{\textbf{Parameter values}} &\multicolumn{2}{c|}{\textbf{Size} ($b=160$)} &\multicolumn{2}{c|}{\textbf{Size} ($b=320$)}&\textbf{Size} \\
$n$ & $\omega$ & $\lambda_1=\lambda_2$ & $\mu_1$ & $\mu_2$ & $\mu^*$  & $T_1$ & $T_2$ & $T_1$ & $T_2$ & $t$ \\ \hline
5000& 0.00247& 2.71& $3.380 \times 10^{-4}$ & $-3.380 \times 10^{-4}$ & 0  & 0.043  &0.043 & 0.046   & 0.054& 0.165 \\ \hline
5000& 0.00247& 2.71& $3.245 \times 10^{-4}$ & $-3.245 \times 10^{-4}$ & 0 & 0.058 &0.045 & 0.057   &0.047 &0.185\\
5000& 0.00247& 2.71&$3.515 \times 10^{-4}$ & $-3.515 \times 10^{-4}$ & 0  &0.043  &0.045 &0.046  &0.048 &0.175 \\
5000& 0.000247& 27.1& $3.380 \times 10^{-4}$ & $-3.380 \times 10^{-4}$ & 0 &0.058 &0.057 &0.056  &0.056 & 0.171 \\
5000& 0.00247& 2.71& 0.3380 & -0.3380 & 0 & 0.056 &0.057 & 0.062  &0.061 &0.222\\
5000& 0.00247& 2.71& 0.3245 & -0.3245 & 0 & 0.055 &0.045 & 0.057  &0.057 &0.232\\
5000& 0.00247& 2.71 & 0.1966 & -0.2023 & 0  &0.059  &0.049   &0.054  &0.049 &0.248\\
5000& 0.00247& 2.71& 0.3515 & -0.3515 & 0&0.057  &0.052 &0.060  &0.064 &0.250\\
5000& 0.00247& 2.71& 0.4793 & -0.4793 & 0 & 0.051 & 0.053  & 0.056 & 0.062 &0.228\\
5000& 0.000247& 27.1& 0.3380 & -0.3380 & 0  &0.054&0.050 &0.059&0.046 &0.233 \\
5000& 0.0000247& 271& 0.3380 & -0.3380 & 0     &0.050  &0.054 &0.050  &0.065 &0.226 \\
10000& 0.00247& 2.71& 0.3245 & -0.3245 & 0 & 0.044  &0.046 & 0.047  &0.050 &0.209\\
10000& 0.00247& 2.71& 0.1966 & -0.1966 & 0  &0.048 &0.050 &0.043 &0.050 &0.213\\
10000& 0.00247& 2.71& 0.3515 & -0.3515 & 0 &0.052 &0.047 &0.060  &0.044 &0.243\\
10000& 0.00247& 2.71& 0.4793 & -0.4793 & 0 & 0.052& 0.048 & 0.055& 0.045&0.235\\\hline
\end{tabular}
\caption{\footnotesize Size of the $t$ and $T$ tests for EACD(1,1) durations with finite variance.} \label{tab:eacdsizeFin}
\end{center}
\end{table}
%\end{landscape}

\begin{landscape}
\begin{table}[H]
\begin{center}
\footnotesize
\begin{tabular}{cccccc|cc|cc}
\hline
\multicolumn{6}{c|}{\textbf{Parameter values}} &\multicolumn{2}{c|}{\textbf{Power} ($b=160$)} &\multicolumn{2}{c}{\textbf{Power} ($b=320$)} \\
$n$ & $\omega$ & $\lambda_1=\lambda_2$ & $\mu_1$ & $\mu_2$ & $\mu^*$  & $T_1$ & $T_2$ & $T_1$ & $T_2$\\ \hline
5000&0.00247&2.71& $3.387 \times 10^{-4}$ & $-3.373 \times 10^{-4}$ & $3.846 \times 10^{-6}$  &0.043 &0.054 &0.052 &0.063 \\
5000&0.00247&2.71& $3.387 \times 10^{-4}$ & $-3.245 \times 10^{-4}$ & $3.846 \times 10^{-5}$  &0.175 &0.196 &0.187 &0.203\\
5000&0.00247&2.71& $3.515 \times 10^{-4}$ & $-3.373 \times 10^{-4}$ & $3.846 \times 10^{-5}$  &0.145 &0.193  &0.155  &0.205 \\
5000&0.000247&27.1& $3.387 \times 10^{-4}$ & $-3.373 \times 10^{-4}$ & $3.846 \times 10^{-5}$ &0.093 &0.096  &0.095  &0.098\\
5000&0.00247&2.71& 0.3387 & -0.3373 & $3.846 \times 10^{-3}$  &0.068 &0.066&0.067 &0.070\\
5000&0.00247&2.71& 0.3387 & -0.3245 & $3.846 \times 10^{-2}$  &0.176 &0.204 &0.178 &0.211 \\
5000&0.00247&2.71& 0.3387 & -0.1966 & $3.846 \times 10^{-1}$   &0.964 &1  &0.954 &0.999  \\
5000&0.00247&2.71& 0.3515 & -0.3373 & $3.846 \times 10^{-2}$  &0.173 &0.219 &0.171 &0.216\\
5000&0.00247&2.71& 0.4793 & -0.3373 & $3.846 \times 10^{-1}$  &0.906 &0.999 &0.918 &0.999 \\
5000&0.000247&27.1& 0.3387 & -0.3373 & $3.846 \times 10^{-2}$  &0.088&0.089 &0.094&0.091\\
5000&0.0000247&271& 0.3387 & -0.3373 & $3.846 \times 10^{-1}$  &0.198 &0.282 &0.190 &0.276 \\
10000&0.00247&2.71& 0.3387 & -0.3245 & $3.846 \times 10^{-2}$  &0.250 &0.311 &0.253 &0.328\\
10000&0.00247&2.71& 0.3387 & -0.1966 & $3.846 \times 10^{-1}$  &0.997 &1  &0.992 &1\\
10000&0.00247&2.71& 0.3515 & -0.3373 & $3.846 \times 10^{-2}$  &0.230 &0.302 &0.243 &0.335\\
10000&0.00247&2.71& 0.4793 & -0.3373 & $3.846 \times 10^{-1}$  &0.972&1 &0.962&1\\ \hline
\end{tabular}
\caption{\footnotesize Power of the $T$ tests for EACD(1,1) durations with finite variance.} \label{tab:eacdpowerFin}
\end{center}
\end{table}
\end{landscape}

%\begin{landscape}
\begin{table}[H]
\begin{center}
\footnotesize
\begin{tabular}{cccccc|ccc}
\hline
\multicolumn{6}{c|}{\textbf{Parameter values}} &\multicolumn{3}{c}{\textbf{Size}}  \\
$n$ & $\omega$ & $\lambda_1=\lambda_2$ & $\mu_1$ & $\mu_2$ & $\mu^*$  & $T_1$ & $T_2$ & $t$ \\ \hline
5000&0.00369&2.71& $1.848 \times 10^{-4}$ & $-1.848 \times 10^{-4}$ & 0  &0.040  &0.041 &0.393 \\
5000&0.00369&2.71& $1.713 \times 10^{-4}$ & $-1.713 \times 10^{-4}$ & 0 &0.043   &0.044 &0.409\\
5000&0.00369&2.71&$1.983 \times 10^{-4}$ & $-1.983 \times 10^{-4}$ & 0  &0.041  &0.040  &0.439 \\
5000&0.000369&27.1& $1.848 \times 10^{-4}$ &$-1.848 \times 10^{-4}$ & 0 &0.042  &0.041  &0.372  \\
5000&0.00369&2.71& 0.1848& -0.1848 & 0 &0.041    &0.042  &0.448\\
5000&0.00369&2.71& 0.1713 & -0.1713 & 0 & 0.043  &0.040  &0.417\\
5000&0.00369&2.71& 0.0435 & -0.0435 & 0  &0.060  &0.054   &0.571\\
5000&0.00369&2.71& 0.1983 & -0.1983 & 0 &0.046  &0.043   &0.448\\
5000&0.00369&2.71& 0.3261 & -0.3261 & 0 & 0.053 & 0.059  &0.506\\
5000&0.000369&27.1& 0.1848& -0.1848 & 0  &0.043&0.041    &0.385 \\
5000&0.0000369&271& 0.1848 & -0.1848 &0     &0.045  &0.045 &0.379 \\
10000&0.00369&2.71& 0.1713 & -0.1713 & 0 & 0.045  &0.040 &0.435\\
10000&0.00369&2.71& 0.0435 & -0.0435 & 0  &0.055 &0.047&0.536\\
10000&0.00369&2.71& 0.1983 & -0.1983 & 0 &0.043  &0.042&0.453\\
10000&0.00369&2.71& 0.3261 & -0.3261 & 0 & 0.046& 0.051&0.517\\\hline
\end{tabular}
\caption{\footnotesize Size of the $t$ and $T$ tests for EACD(1,1) durations with infinite variance.} \label{tab:eacdsizeInf}
\end{center}
\end{table}
%\end{landscape}

%\begin{landscape}
\begin{table}[H]
\begin{center}
\footnotesize
\begin{tabular}{cccccc|cc}
\hline
\multicolumn{6}{c|}{\textbf{Parameter values}} &\multicolumn{2}{c}{\textbf{Power}}  \\
$n$ & $\omega$ & $\lambda_1=\lambda_2$ & $\mu_1$ & $\mu_2$ & $\mu^*$  & $T_1$ & $T_2$  \\ \hline
5000&0.00369&2.71& $1.855 \times 10^{-4}$ & $-1.841 \times 10^{-4}$ & $3.846 \times 10^{-6}$   &0.041 &0.042  \\
5000&0.00369&2.71& $1.855 \times 10^{-4}$ & $-1.713 \times 10^{-4}$ & $3.846 \times 10^{-5}$   &0.060 &0.056 \\
5000&0.00369&2.71&$1.983 \times 10^{-4}$ & $-1.841 \times 10^{-4}$ & $3.846 \times 10^{-5}$   &0.046  &0.047  \\
5000&0.000369&27.1& $1.855 \times 10^{-4}$ & $-1.841 \times 10^{-4}$ & $3.846 \times 10^{-5}$  &0.042  &0.043   \\
5000&0.00369&2.71& 0.1855 & -0.1841 & $3.846 \times 10^{-3}$ &0.046 &0.044\\
5000&0.00369&2.71&0.1855 & -0.1713 & $3.846 \times 10^{-2}$  &0.047 &0.046\\
5000&0.00369&2.71&0.1855& -0.0435 & $3.846 \times 10^{-1}$   &0.408 &0.622    \\
5000&0.00369&2.71&0.1983 & -0.1841 & $3.846 \times 10^{-2}$ &0.056 &0.058\\
5000&0.00369&2.71&0.3261 & -0.1841 & $3.846 \times 10^{-1}$ &0.244 &0.299 \\
5000&0.000369&27.1&0.1855 & -0.1841 & $3.846 \times 10^{-2}$  &0.044&0.045 \\
5000&0.0000369&271&0.1855 & -0.1841 & $3.846 \times 10^{-1}$   &0.057 &0.059 \\
10000&0.00369&2.71&0.1855 & -0.1713 & $3.846 \times 10^{-2}$  &0.043 &0.046\\
10000&0.00369&2.71& 0.1855 & -0.0435 & $3.846 \times 10^{-1}$  &0.453 &0.682\\
10000&0.00369&2.71& 0.1983 & -0.1841 & $3.846 \times 10^{-2}$  &0.050 &0.052\\
10000&0.00369&2.71& 0.3261 & -0.1841 & $3.846 \times 10^{-1}$ &0.258&0.358\\\hline
\end{tabular}
\caption{\footnotesize Power of the $T$ tests for EACD(1,1) durations with infinite variance.} \label{tab:eacdpowerInf}
\end{center}
\end{table}
%\end{landscape}

\subsection{LMSD durations}
\label{subsec:lmsdsimu}

We generate $\{\tau_{k}^{(i)}\}$ from the LMSD model
\[ \tau_{k}=\epsilon_{k} e^{Y_{k}} \] where $\{\epsilon_{k}\} \iid Weibull(\delta,\gamma)$, and
($\delta$, $\gamma$) are the scale and shape parameter; $\{Y_{k}\}$ is a Gaussian
\textit{ARFIMA}$(1,d,0)$ process with innovations $w_{k} \iid N(0, \sigma_w ^{2})$, given by
\[ (1-\alpha L)(1-L)^{d}Y_{k}=w_{k} \] with $|\alpha|<1$, and $d$ is the long memory parameter with
$d =H-1/2 \in (0,1/2)$. We generate models without and with microstructure shocks
$\{\eta_{k}\}$. For all of these models, we find that $t$ is oversized and therefore we will not
examine its power. We will also study the size and power of the $T$ tests and the leverage effect
induced by the microstructure shocks. In all simulations, we use 1000 realizations and block size
$b=160$.

\subsubsection{Parameter calibration}
\label{subsec:lmsdParCali}

For both LMSD models without and with $\{\eta_{k}\}$, we consider a total of 30 configurations of
parameters and sample size, 15 for evaluating size (Table \ref{tab:lmsdsize}), and 15 for evaluating
power (Table \ref{tab:lmsdpower}). Parameters for the LMSD model are ($\alpha$, $\delta$, $\gamma$,
$\sigma_w ^{2}$, $d$, $\sigma$, $\mu_1$, $\mu_2$). Note that we assume that ($\alpha$, $\delta$,
$\gamma$, $\sigma_w ^{2}$, $d$, $\sigma$) are the same for both processes $\{\tau_{k}^{(1)}\}$,
$\{\tau_{k}^{(2)}\}$.

\noindent \textit{LMSD without microstructure shocks}

\begin{figure}[H]
\begin{center}
\includegraphics[scale=.8, trim = 35mm 140mm 30mm 20mm]{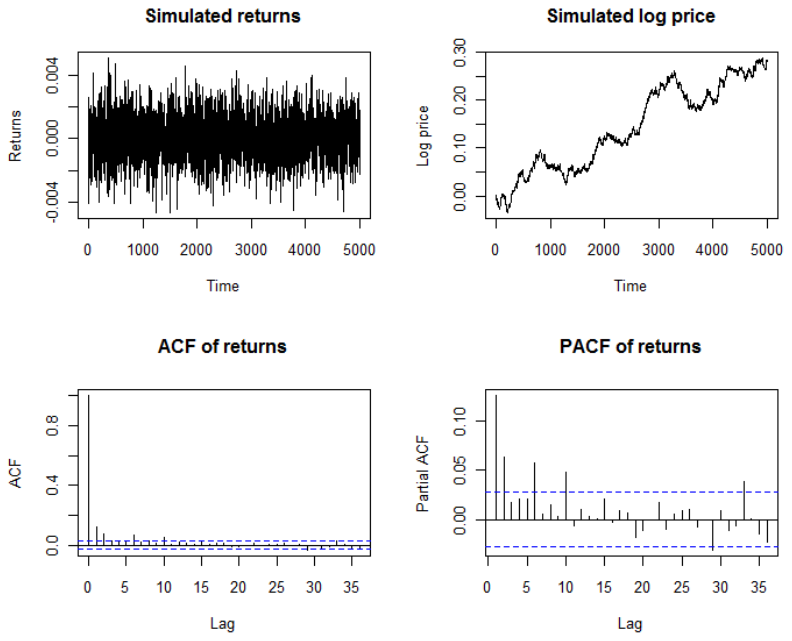}
  \caption{\footnotesize One simulated realization of the calibrated LMSD model without
    $\{\eta_{k}\}$ (first row in Table \ref{tab:lmsdpower}). Time series plot of returns and log
    price, as well as ACF and PACF, $n$=5000. The model is calibrated using Boeing's and Fama/French
    factor Rm-Rf, such that $\alpha=-.4212$, $\delta=0.2$, $\gamma=1.3376$, $\sigma_w ^{2}=0.2368$,
    $d=0.3545$, $\sigma=0.0003$, and the mean and standard deviation for the 5-minute returns are
    $\mu^*= 3.846\times 10^{-6}$, $\tilde{\sigma} = 0.0012$. }\label{fig:simuLMSD}
\end{center}
\end{figure}

 For the first configuration corresponding to power (Table \ref{tab:lmsdpower}), we calibrate the model to match the ($\alpha, \delta, \gamma, \sigma_w ^{2}, d$) observed in the Boeing series used in \citet{deo:hsieh:hurvich:2010}. We use the same $\sigma$ as for the exponential durations and set $\mu_1$ and $\mu_2$ to match the mean and standard deviation of the Rm-Rf. This calibrated model also accounts for the discrepancies in the time units, and is given in the first row of Table \ref{tab:lmsdpower}. Similarly as for the EACD(1,1) model, for the LMSD model, the functional relationship between the model parameters and the mean and variance of the calendar-time returns also depend on the variance of counts. We estimate the variance of counts using simulation and set $\mu_1$ and $\mu_2$ accordingly. Figure~\ref{fig:simuLMSD}, for one simulated realization from this model with $n=5000$, shows the time series plots of returns and log price, as well as the ACF and PACF of the returns.

For the remaining configurations in Table \ref{tab:lmsdpower}, we increase $\mu^*$ by varying $\delta$ (which changes $\lambda$) and $\mu_i$. We consider two sample sizes, $n=5000$ and $n=10000$. The power should generally increase as $\mu^*$ increases.

For the first configuration corresponding to size (Table \ref{tab:lmsdsize}), we set the mean to zero. For the remaining configurations in these two tables, we vary the $\delta$ and $\mu_i$ similarly as in Table \ref{tab:lmsdpower} while keeping $\mu^*=0$. \\

\noindent \textit{LMSD with microstructure shocks}

For $i=1,2$, define the microstructure shock $\eta_k$ as
\[
\eta_k = Y_{k-1}-Y_k \;.
\]

With $\{\eta_{k}\}$, the functional relationship between the model parameters and the variance of
the calendar-time returns now depend on not only the variance of counts, but also the variance of
$\{\eta_{k}\}$ and the covariance of counts and $\{\eta_{k}\}$. Since it is difficult to calibrate
this variance and covariance, for simplicity, we will use the same values of all parameters
($\alpha$, $\delta$, $\gamma$, $\sigma_w ^{2}$, $d$, $\sigma$, $\mu_1$, $\mu_2$) as those used in
the LMSD models without $\{\eta_{k}\}$. The size and power results are given in Tables
\ref{tab:lmsdsize} (size) and \ref{tab:lmsdpower} (power).

% We use the same values of ($\alpha, \gamma, \sigma_w ^{2}, d, \sigma$) as for all LMSD models
% without $\{\eta_{k}\}$. To set $\mu_1$ and $\mu_2$, note that now the functional relationship
% between the model parameters and the variance of the calendar-time returns depend on not only the
% variance of counts, but also the variance of $\{\eta_{k}\}$ and the covariance of counts and
% $\{\eta_{k}\}$. Since it is difficult to calibrate the variance and covariance related to the
% microstructure, for simplicity, we would also use the same values of $\mu_1$ and $\mu_2$ as those
% of all LMSD models without $\{\eta_{k}\}$. The size and power results are also given in Tables
% \ref{tab:lmsdsize} (size) and \ref{tab:lmsdpower} (power).

\begin{figure}[H]
\begin{center}
  \includegraphics[scale=.8, trim = 30mm 142mm 30mm 50mm]{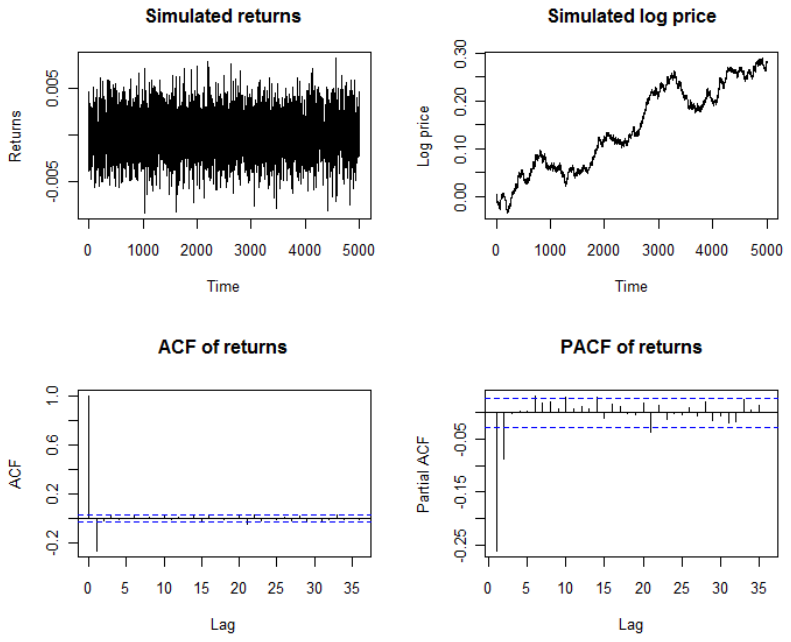}
  \caption{\footnotesize One simulated realization of the calibrated LMSD model with $\{\eta_{k}\}$
    (first row in Table \ref{tab:lmsdpower}). Time series plot of returns and log price, as well as
    ACF and PACF, $n$=5000. The model is calibrated using Boeing's and Fama/French factor Rm-Rf,
    such that $\alpha=-.4212$, $\delta=0.2$, $\gamma=1.3376$, $\sigma_w ^{2}=0.2368$, $d=0.3545$,
    $\sigma=0.0003$, and the mean and standard deviation for the 5-minute returns are $\mu^*=
    3.846\times 10^{-6}$, $\tilde{\sigma} = 0.0012$. }\label{fig:simuLMSDeta}
\end{center}
\end{figure}

\begin{figure}[H]
\begin{center}
  \includegraphics[scale=.7, trim = 30mm 121mm 30mm 44mm]{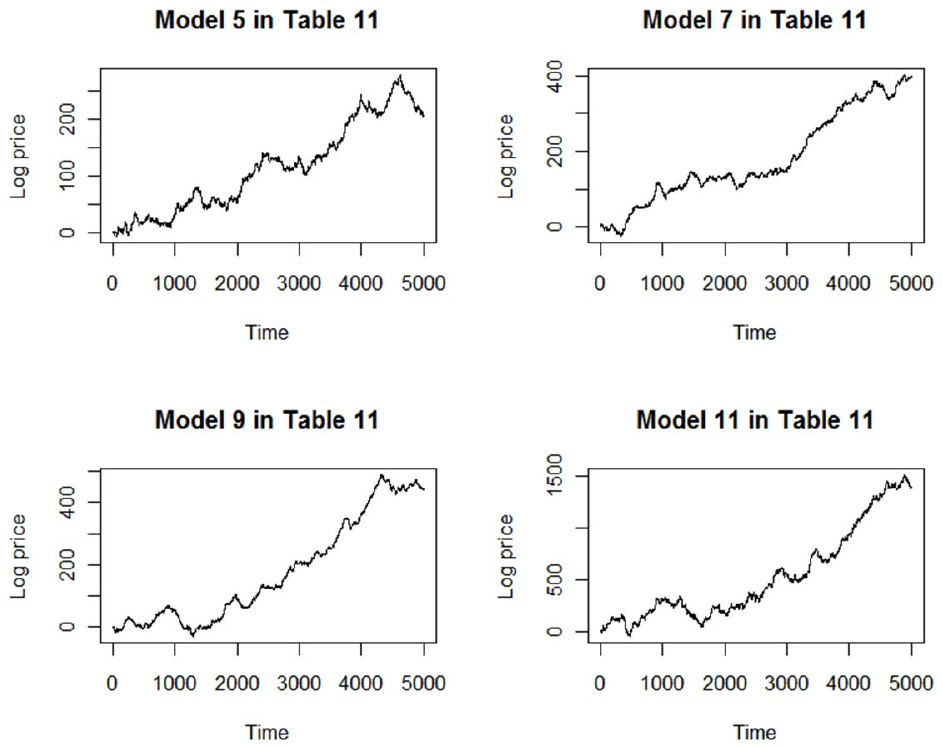}
  \caption{\footnotesize Simulated log price for models with $\{\eta_{k}\}$ in rows 5 ($\mu^*=3.846
    \times 10^{-3}$), 7, 9, and 11 (all with $\mu^*=3.846 \times 10^{-1}$) of Table
    \ref{tab:lmsdpower}.}\label{fig:simuLMSDeta2}
\end{center}
\end{figure}

Figure~\ref{fig:simuLMSDeta}, for one simulated realization from first configuration of Table
\ref{tab:lmsdpower} with $\{\eta_{k}\}$, shows the time series plots of returns and log price, as
well as the ACF and PACF of the returns. Figure \ref{fig:simuLMSDeta2} shows the time series plots
of the log price for the models with $\{\eta_{k}\}$ of rows 5 ($\mu^*=3.846 \times 10^{-3}$), 7, 9,
and 11 (all with $\mu^*=3.846 \times 10^{-1}$) in Table \ref{tab:lmsdpower}.

\subsubsection{Leverage effect}
\label{subsec:leverage}

In model (\ref{eq:UnivariateModel}), we assume that the microstructure shocks $\{\eta_{i,k}\}$ are
independent of the efficient shocks $\{e_{1,k}\}$, $\{e_{2,k}\}$, but not necessarily of the
counting process $(N_1,N_2)$. This allows for a leverage effect. See
\citet{aue:hurvich:horvath:soulier:2014}.  The leverage is then defined as the correlation between
current calendar-time return $r_k$ and absolute value of the next calendar-time return $|r_{k+1}|$
\[ leverage = corr(r_k, |r_{k+1}|) \;.\]

We test the model with different memory parameter $d=.23$ and $d=.4$, keeping other parameters the
same as in the first configuration of Table \ref{tab:lmsdpower}. We generate 100 realizations, with
$n =5000$, calculate sample leverage for each realization, and use the $t$-test for the mean,
\begin{align*}
  &H_0: \; leverage = 0 \\
  &H_1: \; leverage < 0 \;.
\end{align*}

Table \ref{tab:leverage} shows the mean of sample leverage and the $p$-value for models with
different memory parameters $d$. Without the $\{\eta_{k}\}$, there is no leverage in returns, while
including the $\{\eta_{k}\}$ can induce a leverage effect, i.e. negative correlation between current
return and the next period's absolute return.

\begin{table}[H]
\begin{center}
\footnotesize
\begin{tabular}{c|cc|cc}
  \hline
  \textbf{Model} & \multicolumn{2}{c|}{$d=0.23$} &\multicolumn{2}{c}{$d=0.4$}  \\
  &\textbf{Mean} & \textbf{p-value}&\textbf{Mean} & \textbf{p-value} \\ \hline
  \textbf{Returns without $\eta_k$}& 0.00062 & 0.6265 & -0.00374& 0.3727 \\
  \textbf{Returns with $\eta_k$}&-0.02299& 0.0000 &-0.01883 & 0.0000 \\ \hline
\end{tabular}
\caption{\footnotesize Means and p-values of sample leverage for models with $d=.23$ and
  $d=.4$.} \label{tab:leverage}
\end{center}
\end{table}

\subsubsection{Size/Power of the $t$ and $T$ tests for LMSD durations}
\label{subsec:lmsdSizePower}

Tables \ref{tab:lmsdsize} and \ref{tab:lmsdpower} show that there is no significant difference in
the size and power between models without and with microstructure shocks $\{\eta_{k}\}$. For both
models, Table \ref{tab:lmsdsize} shows that the $t$ test is oversized. $T_1$ is correctly sized in
all cases, while $T_2$ is a bit oversized. As we increase $\lambda$ or $n$, the size of $T_2$
becomes smaller. Table \ref{tab:lmsdpower} shows that the power of the $T$ test performs similarly
as in the Poisson model: the power increases as $\mu^*$ increases; as $n$ is increased holding all
model parameters fixed, the power increases; $T_2$ has higher power than $T_1$.

\begin{landscape}
\begin{table}[H]
\begin{center}
\footnotesize
\begin{tabular}{cccccc|ccc|ccc}
\hline
\multicolumn{6}{c|}{\textbf{Parameter values}} &\multicolumn{3}{c|}{\textbf{Size} (without $\eta_{k}$)} &\multicolumn{3}{c}{\textbf{Size} (with $\eta_{k}$)}  \\
$n$ & $\delta$ & $\lambda_1=\lambda_2$ & $\mu_1$ & $\mu_2$ & $\mu^*$  & $T_1$ & $T_2$ & $t$ & $T_1$ & $T_2$ & $t$ \\ \hline
5000&0.2 &4.58& $2.293 \times 10^{-4}$ & $-2.293 \times 10^{-4}$ & 0  & 0.057  &0.071 &0.307 & 0.051  &0.080&0.307  \\
5000&0.2& 4.58& $2.213 \times 10^{-4}$ & $-2.213 \times 10^{-4}$ & 0 & 0.059 &0.080 &0.312 & 0.057 &0.084 &0.311\\
5000&0.2& 4.58& $2.373 \times 10^{-4}$ & $-2.373 \times 10^{-4}$ & 0  &0.052  &0.064   &0.291 &0.061  &0.079   &0.320 \\
5000&0.02&45.8& $2.293 \times 10^{-4}$ & $-2.293 \times 10^{-4}$ & 0 &0.052  &0.083  &0.337 &0.059  &0.072 &0.314  \\
5000&0.2& 4.58& 0.2293 & -0.2293 & 0 &0.058   &0.073&0.341 &0.061  &0.063&0.345\\
5000&0.2& 4.58& 0.2213 & -0.2213 & 0 & 0.063  &0.088&0.352 & 0.057  &0.089&0.351\\
5000&0.2& 4.58 &0.1458& -0.1458 & 0  &0.046  &0.073    &0.376 &0.044  &0.081   &0.341\\
5000&0.2& 4.58&0.2373 & -0.2373 & 0 &0.057  &0.078&0.366 &0.053  &0.076&0.339\\
5000&0.2& 4.58&0.3129 & -0.3129 & 0 & 0.056 & 0.077  &0.343 & 0.062 & 0.077 &0.372\\
5000&0.02& 45.8& 0.2293 & -0.2293 & 0 &0.052&0.070&0.357 &0.054&0.071&0.364 \\
5000&0.002& 458 &0.2293 & -0.2293 &0     &0.054  &0.057&0.371 &0.053  &0.052&0.389 \\
 10000&0.2& 4.58&0.2213 & -0.2213 & 0 & 0.057  &0.068&0.366 & 0.062  &0.070&0.368\\
 10000&0.2& 4.58& 0.1458 & -0.1458  & 0  &0.053 &0.070&0.368 &0.056 &0.071&0.370\\
10000&0.2& 4.58& 0.2373 & -0.2373 & 0 &0.048  &0.071&0.368 &0.051  &0.070&0.376\\
10000&0.2& 4.58& 0.3129 & -0.3129 & 0 & 0.052& 0.071&0.373 & 0.051& 0.071&0.374\\\hline
\end{tabular}
\caption{\footnotesize Size of the $t$ and $T$ tests for LMSD durations.} \label{tab:lmsdsize}
\end{center}
\end{table}
\end{landscape}

\begin{landscape}
\begin{table}[H]
\begin{center}
\footnotesize
\begin{tabular}{cccccc|cc|cc}
\hline
\multicolumn{6}{c|}{\textbf{Parameter values}} &\multicolumn{2}{c|}{\textbf{Power} (without $\eta_{k}$)} &\multicolumn{2}{c}{\textbf{Power} (with $\eta_{k}$)}  \\
$n$ & $\delta$ & $\lambda_1=\lambda_2$ & $\mu_1$ & $\mu_2$ & $\mu^*$  & $T_1$ & $T_2$ & $T_1$ & $T_2$ \\ \hline
5000&0.2 &4.58& $2.297 \times 10^{-4}$ & $-2.289 \times 10^{-4}$ & $3.846 \times 10^{-6}$ &0.072 &0.080&0.071&0.090 \\
5000&0.2& 4.58& $2.297 \times 10^{-4}$ & $-2.213 \times 10^{-4}$ & $3.846 \times 10^{-5}$ & 0.134&0.202 &0.128&0.197\\
5000&0.2& 4.58& $2.373 \times 10^{-4}$ & $-2.289 \times 10^{-4}$ & $3.846 \times 10^{-5}$ &0.116&0.161 &0.130&0.204 \\
5000&0.02&45.8& $2.297 \times 10^{-4}$ & $-2.289 \times 10^{-4}$ & $3.846 \times 10^{-5}$ &0.076&0.092 &0.075 &0.087\\
5000&0.2& 4.58& 0.2297 & -0.2289 & $3.846 \times 10^{-3}$   &0.074 &0.08 &0.067&0.069\\
5000&0.2& 4.58&0.2297 & -0.2213 & $3.846 \times 10^{-2}$  &0.139 &0.193 &0.120 &0.192\\
5000&0.2& 4.58 &0.2297& -0.1458 & $3.846 \times 10^{-1}$   &0.881 &0.999 &0.905 &0.998    \\
5000&0.2& 4.58&0.2373 & -0.2289  & $3.846 \times 10^{-2}$  &0.141 &0.193 &0.107 &0.166\\
5000&0.2& 4.58&0.3129 & -0.2289  & $3.846 \times 10^{-1}$ &0.697 &0.965 &0.688 &0.974 \\
5000&0.02& 45.8& 0.2297& -0.2289 & $3.846 \times 10^{-2}$  &0.079&0.089  &0.062&0.080 \\
5000&0.002& 458 &0.2297 & -0.2289 &$3.846 \times 10^{-1}$      &0.091 &0.105 &0.083 &0.098 \\
 10000&0.2& 4.58&0.2297 & -0.2213 & $3.846 \times 10^{-2}$  &0.154 &0.209 &0.144 &0.204\\
 10000&0.2& 4.58& 0.2297 & -0.1458 & $3.846 \times 10^{-1}$  &0.940 &1 &0.949 &1\\
10000&0.2& 4.58& 0.2373 & -0.2289 & $3.846 \times 10^{-2}$  &0.145 &0.205 &0.143 &0.195\\
10000&0.2& 4.58& 0.3129 & -0.2289 & $3.846 \times 10^{-1}$ &0.793&0.997 &0.780&0.990 \\\hline
\end{tabular}
\caption{\footnotesize Power of the $T$ tests for LMSD durations.} \label{tab:lmsdpower}
\end{center}
\end{table}
\end{landscape}

\section{Data analysis}
\label{sec:data}

We study the daily returns of the Fama/French factor, Rm-Rf, from Kenneth French's data
library. According to the website, ``Rm-Rf is the excess return on the market, value-weight return
of all CRSP (Center for Research in Security Prices) firms incorporated in the US and listed on the
NYSE, AMEX, or NASDAQ that have a CRSP share code of 10 or 11 at the beginning of month t, good
shares and price data at the beginning of t, and good return data for t minus the one-month Treasury
bill rate (from Ibbotson
Associates)"(\url{http://mba.tuck.dartmouth.edu/pages/faculty/ken.french/data_library.html}). The
data ranges from July 1, 1926 to December 31, 2013, a total 23133 daily excess returns with sample
mean 0.000285.

\begin{figure}[H]
  \begin{center}
    \includegraphics[scale=.8, trim = 30mm 130mm 30mm 40mm]{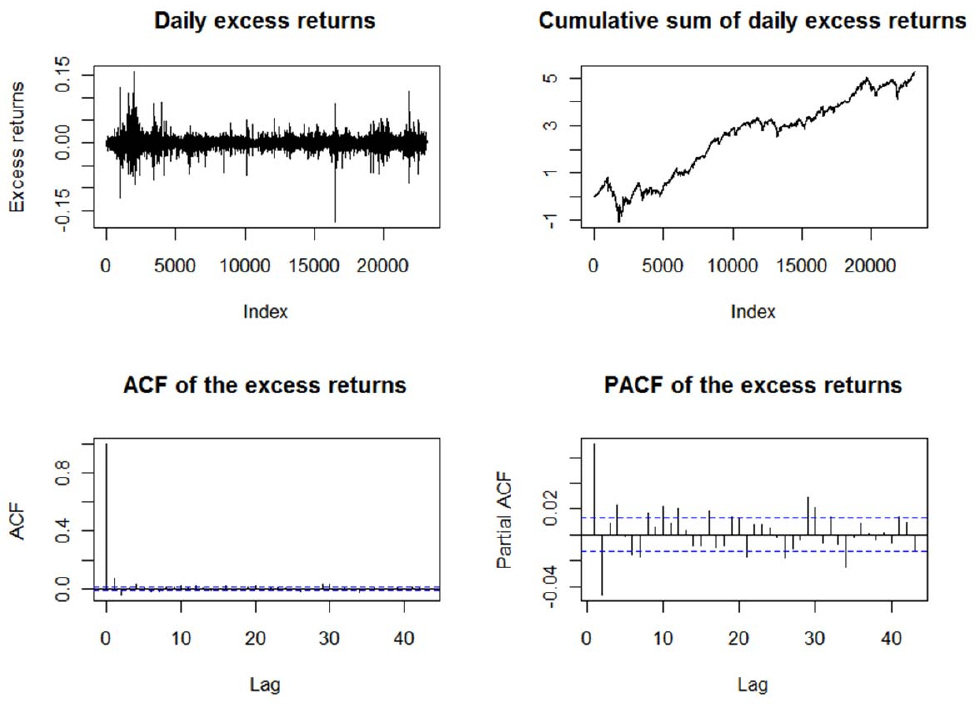}
    \caption{\footnotesize Time series plot, log price, ACF, and PACF of daily excess returns from
      July 1, 1926 to December 31, 2013.}\label{fig:dailyexcessreturn}
\end{center}
\end{figure}

Figure~\ref{fig:dailyexcessreturn} shows the time series plot, log price, ACF and PACF of daily
excess returns. Both ACF and PACF show statistically significant lags. We perform a $t$-test based
on the usual standard error and the Newey-West standard error, which allows for serial correlation,
and also our $T_1$ test for the following hypothesis:
\begin{align*}
&H_0: \mu^* = 0\\
&H_1: \mu^* > 0
\end{align*}
where $\mu^*$ is the expected excess return (equity premium).  To conduct the $T_1$ test, we
calculate the $T_1$ statistic for the entire data set, and use the same subsampling procedure as
defined in Section~\ref{sec:simulation} to obtain the empirical quantiles of the limiting
distribution of T. As suggested by \citet{jach:mcelroy:politis:2012}, with the sample size of
$n=23133$, we consider the block size $b\in \{732, 976, 1302, 1736, 2314, 3086, 4116, 5488\}$.  This
covers subsample sizes ranging from $3-24\%$ of the sample size. The results are given in Table
\ref{tab:Ttest}.

At the 0.05 significant level, both the ordinary $t$-test and the Newey-West $t$-test reject the
null hypothesis with very small p-values, 2.59e-05 and 7.40e-05 respectively. This seems to provide
extremely strong evidence for an equity premium. However, the much larger $p$-values for all of the
$T_1$ tests considered indicate far weaker evidence against the null hypothesis. This finding is
consistent with the simulation results in Section \ref{sec:simulation} for the EACD and LMSD
durations. The $p$-values corresponding to the ordinary $t$-tests may be spuriously low, and the
ones corresponding to the $T_1$ tests may be more reliable.

\begin{table}[H]
\begin{center}
\footnotesize
\begin{tabular}{ccccc}
\hline
	\multicolumn{5}{c}{\textbf{T test}} \\ \hline
\textbf{Block size (b)}& 	732	&	977	&	1302	&	1736	 \\
\textbf{$p$-value}&	0.054	&	0.057	&	0.056	&	0.058\\ \hline
\textbf{Block size (b)}&	2351	&	3087  	&	4117	&	5489	\\
\textbf{$p$-value}	&	0.064	&	0.062  	&	0.065	&	0.057 \\\hline
\end{tabular}
\caption{\footnotesize $T_1$ test results for daily excess returns Rm-Rf.} \label{tab:Ttest}
\end{center}
\end{table}

\section{Discussion: Long Memory and Heavy Tails of Stock Returns}
\label{sec:discussion}

The introduction of a nonzero mean in the efficient shocks in the model
(\ref{eq:UnivariateModel}) provides a link by which properties of intertrade
durations can affect those of certain quantities that are observed at a
macroscopic level. We have focused so far on inference for the trend (based on
studying the asymptotic distribution of the log price). To illustrate just one
of the variety of possible additional quantities of interest, we now turn our
attention to properties of returns.

\citet{lo:1991} investigated whether stock returns have long memory, and
\citet{Mandelbrot:1963} argued that returns have infinite variance. Both of
these propositions have met with considerable controversy, but under the model
(\ref{eq:trend}) both could contain an important grain of truth. Generalizing
the analysis presented so far leads to a more nuanced interpretation of what
these propositions could mean.

From here on in this section, when we mention sequences of random variables, we
allow for suitable renormalization (centering and scaling) without always
specifically mentioning or writing the renormalization. So the discussion here
is somewhat informal, but can be made mathematically rigorous. We focus here on
the case $\gamma >1/2$.

Proposition \ref{prop:weakconv-returns} implies that partial sums of returns
(after suitable renormalization) converge in distribution to a random variable
that need not be Gaussian. This theorem has allowed us to discuss issues related
to inference for the slope parameter, which is the expectation of the average
return. It is also of interest to ask if one can go further and say something
about the joint distribution of the returns themselves, rather than their sum. Therefore, we will now discuss the joint
limiting distribution of any fixed number of contiguous returns at long
horizons.

Although we have so far taken the time spacing in defining the returns to be 1,
there is no essential reason for this and here we replace it by an arbitrary $T
>0$, and we define the returns with respect to this time spacing as $r_{j,T} =
y(jT)-y((j-1)T)$. Now consider the first $M$ of these returns, where $M$ is
fixed. It follows from our assumptions here (by arguments similar to the proof
of Proposition \ref{prop:weakconv-returns} and by Theorem 7.3.2 of
\citet{whitt:2002}) that the joint distribution of these $M$ returns (after
suitable renormalization) converges as $T\rightarrow \infty$ to the distribution
of $M$ contiguous increments of the limiting process.

In our LMSD example, assuming finite variance and an exponential volatility
function, the limiting process is fractional Brownian motion. Thus in this case,
the $M$ returns converge in distribution to $M$ contiguous observations of a
fractional Gaussian noise. In this sense, it could be said that the returns
(computed at a sufficiently high level of aggregation) have long
memory. Simulations not shown here of the model (\ref{eq:trend}) in this LMSD
case show that it may be hard to detect this long memory due to the additive
noise that arises from the second term on the righthand side of
(\ref{eq:trend}).

In the ACD example (and certain cases of the LMSD example as well), it turns out
that the limiting process can be a stable process, for which the increments are
independent and have infinite-variance stable distributions. So here, our $M$
long-horizon returns converge in distribution (as $T \rightarrow \infty$, and
after suitable renormalization) to a sequence of $M$ i.i.d. stable random
variables.  This would seem to correspond to the proposition that returns have
infinite variance. But actually, the truth here may be more subtle.  It can
happen that, for each fixed $T$ the variance of the returns is finite. See, for
example, \citet{whitt:2002} for the underlying point process theory under
heavy-tailed durations. It is even possible to construct an example where
durations have finite variance and still the limit of partial sums of durations
is a stable process, so the returns would once again have finite variance but
converge in distribution to i.i.d. stable random variables with infinite
variance. Such an example may come from durations that obey a positive version
of the renewal-reward process discussed in \citet{taqqu:levy:1986} (see also
\citet{hsieh:hurvich:soulier:2007}). In such a model, durations would have
finite variance but their sums would converge to a process with infinite
variance.

In the case where the limiting process is a L\'evy-stable process, it is of
interest to note that such continuous-time processes have discontinuities with
probability 1. These may correspond to what practitioners refer to as "jumps" in
the log price process, even though under our model the log price process is a
pure-jump process so that all activity consists of jumps.

The main message of this paper is that even in the simple transaction-level
model (\ref{eq:UnivariateModel}) there is a wide variety of possible behaviors
of macroscopic quantities of interest.  Some additional quantities we hope to
study in future work based on this and similar models include: regression
coefficients as used in the market model, estimated cointegrating parameters
(which were considered without a trend term in
\citet{aue:hurvich:horvath:soulier:2014}), and sample autocorrelations.

\section{Proofs}
\label{sec:math}

\begin{proof}[Proof of Proposition~\ref{prop:weakconv-returns}]
  We start by proving the convergences~(\ref{eq:conweak-gamma>1/2})
  and~(\ref{eq:conweak-1/2}). Defining $\zeta_{i,k} =
  e_{i,k}+\eta_{i,k}$, we can write
  \begin{align*}
    y(n) - \mu^* n = \mu_1 \{N_1(n) - \lambda_1 n\} + \mu_2 \{N_2(n) - \lambda_2
    n\} + \sum_{k=1}^{N_1(n)} \zeta_{1,k} + \sum_{k=1}^{N_2(n)} \zeta_{2,k}
  \end{align*}
  Assumption~\ref{hypo:weakconv-counts} implies that $N_i(n)/n \plim \lambda_i$.
  For $t\geq0$, define $x_{i,n}(t) = n^{-1/2} \sum_{k=1}^{N_i(nt)} \zeta_{i,k}$.
  The present assumptions imply assumptions 2.1, 2.2, 2.3 and 2.4 of
  \citet{aue:hurvich:horvath:soulier:2014}, thus Theorem~3.1 therein implies that
0  $x_{i,n} \fidi \sqrt \lambda \sigma_i B_i$, where $B_i$ are mutually
  independent standard Brownian motions.  Since the sequences $\{e_{i,k}\}$ are
  independent of the point processes $N_i$, the sequences of processes $x_{i,n}$
  and $\{n^{-\gamma} (N_i(nt)-n\lambda_i t),t\geq0\}$, $i=1,2$ converge
  jointly. Thus, if $\gamma>1/2$, (\ref{eq:conweak-gamma>1/2}) holds and if
  $\gamma=1/2$, (\ref{eq:conweak-1/2}) holds with $\sigma^2 = \lambda_1
  \sigma_1^2+\lambda_2 \sigma_2^2$.
\end{proof}

\begin{proof}[Proof of Theorem~\ref{theo:ratio}]
  The convergence~(\ref{eq:conv-newstat}) is a straightforward consequence
  of Proposition~\ref{prop:weakconv-returns} and the continuous mapping
  Theorem. If $\mu^*\ne\mu^*_0$, write
  \begin{align*}
    \newstat_n = \frac{\bar{r}(n) - \mu^*} {|\bar{r}(n)-\bar{r}(n/2)|} +
    \frac{\mu^* - \mu_0^*} {|\bar{r}(n)-\bar{r}(n/2)|} \; .
  \end{align*}
  The first term converges weakly to $\limproc(1)/|\limproc(1)-2\limproc(1/2)|$
  and the second term converges in probability to $+\infty$ if $\mu^*>\mu_0^*$
  and to $-\infty$ if $\mu^*<\mu_0^*$.
\end{proof}
\begin{proof}[Proof of Lemma~\ref{lem:consistency-s}]
  \begin{align}
  \label{eq:doubleproduct}
    s_n^2 = \frac1{n-1} \sum_{j=1}^n (r_j-\mu^*)^2 + \frac n{n-1}
    (\bar r_n-\mu^*)^2 + 2 (\mu^* - \bar r) \sum_{j=1}^n (r_j - \mu^*)
\end{align}
Under Assumption~\ref{hypo:weakconv-counts}, the
convergences~(\ref{eq:conweak-gamma>1/2}) and~(\ref{eq:conweak-1/2})
imply that the second term in~(\ref{eq:doubleproduct}) is $o_P(1)$.
Since
\begin{align*}
  r_j - \mu^* & = \sum_{i=1}^2 \sum_{k=N_i(j-1)+1}^{N_i(j)}
  (\mu_i+\zeta_{i,k}) -\lambda_i\mu_i \; ,
\end{align*}
by Assumption~(\ref{eq:condition-consistency-s}), ergodicity of the
marked point processes and mutual independence, we have
\begin{align*}
  \frac 1n \sum_{j=1}^n (r_j-\lambda\mu)^2 \plim \sum_{i=1}^2 \esp
  \left[ \left(\sum_{k=1}^{N_i(1)} (\mu_i+\zeta_{i,k}) -\lambda_i\mu_i
    \right)^2 \right] \; .
 \end{align*}
 Since the last term in~(\ref{eq:doubleproduct}) is a cross product,
 it is therefore also $o_P(1)$.
\end{proof}

\begin{lemma}
  \label{lem:conv-cox}
  Let $\xi$ be a stationary random measure such that $n^{-\gamma} \{\xi(nt) -
  \lambda nt\} \fidi \limproc(t)$.  Let $N$ be a Cox process with stochastic
  intensity $\xi$. Then $n^{-\gamma} \{N(nt)-\lambda nt\} \fidi \limproc(t)$.
\end{lemma}

\begin{proof}
  We only prove the convergence in distribution of $n^{-\gamma}\{\xi(n) -\lambda
  n\}$ to $\limproc(1)$. The convergence of the finite dimensional
  distribution is proved similarly. By conditioning on~$\xi$, we have,
  for all $z\in\Rset$,
  \begin{align*}
    \esp \left[ \rme^{\rmi z t^{-\gamma} \{N(t)-\lambda t\}} \right] & = \esp
    \left[ \rme^{\xi(t) \{\rme^{\rmi z t^{-\gamma}}-1\} - \rmi z t^{1-\gamma} }  \right]  \\
    & = \esp \left[ \rme^{\rmi z t^{-\gamma} \{\xi(t) - \lambda t\}} \right] +
    \esp \left[ \rme^{\rmi z t^{-\gamma} \{\xi(t) - \lambda t\}} \left\{
        \rme^{\xi(t) \{\rme^{\rmi z t^{-\gamma}}-1-\rmi z t^{-\gamma} \} } - 1 \right\} \right] \\
    & = \phi_t(z) + \reste_t(z) \; .
  \end{align*}
  By assumption, $\lim_{t\to\infty} \phi_t(z)=\esp[\rme^{\rmi z
    \limproc(1)}]$. Then, denoting $h(t,z) = \rme^{\rmi z t^{-\gamma}}-1-\rmi z
  t^{-\gamma} $, the remainder term is bounded by
\begin{align*}
  |\reste_t(z)| \leq \esp \left[ \left| \rme^{h(t,z) \xi(t)} - 1 \right| \right]  \; .
\end{align*}
For each $z$, $h(t,z) = O(t^{-2\gamma})$, so the assumption on $\xi$ implies that
$\rme^{h(t,z) \xi(t)} - 1 \plim 0$. Moreover, for all $t>0$ and $z\in\Rset$,
$| \rme^{h(t,z)\xi(t)} - 1 | \leq 1 + \rme^{(cos(z)-1) \xi(t)} \leq 2$.
Thus, the bounded convergence theorem yields that \mbox{$\lim_{t\to\infty}
  \esp[|\rme^{h(t,z)\xi(t)} - 1|] = 0$}.
\end{proof}

We recall the definition of $\theta$-weak dependence from \citet{bardet:doukhan:lang:ragache:2008}.
% \citet{dedecker:etal:2007} and \citet{doukhan:louhichi:1999} For each positive
% integer~$u$, let $\mathcal{F}_u$ be the class of bounded measurable functions
% $f: \Rset^u \rightarrow \Rset$ equipped with the supremum norm
% $\|\cdot\|_\infty$.
For each positive integer~$v$, let $\Rset^v$ be equipped with the $l^1$-norm
\[
\|\mathbf x - \mathbf y \|_1 = \sum\limits_{i=1}^{v} |x_i - y_i|, \; \mathbf x,
\mathbf y \in \Rset^v \;.
\]
Let $\mathcal{G}_v$ be the class of bounded functions $g: \Rset^v \rightarrow
\Rset$ such that
\[
\lip(g) = \sup_{\mathbf x \neq \mathbf y} \frac{|g(\mathbf x)-g(\mathbf y)|}{
  \|\mathbf x-\mathbf y\|_1} < \infty \;.
\]
The quantity in the left hand side is denoted $\lip(g)$ and called the Lipschitz
modulus of $g$. A stationary sequence $\{X_k,k\in\Zset\}$ is said to be
$\theta$ weakly dependent with rate $\{\theta_k,k\geq0\}$ if for any $u, v \in
\mathbb{N}^*$, $(i_1, \dots, i_u) \in \Zset$, $(j_1, \dots, j_v) \in
\Zset$ with $i_1 < \dots< i_u \leq i_u+k \leq j_1 < \cdots <j_v$, and $f \in
\mathcal{G}_u$, $g \in \mathcal{G}_v$,
\begin{align*}
  \cov(f(r_{i_1}, \dots, r_{i_u}), g(r_{j_1}, \dots, r_{j_v})) \leq v \,
  \|f\|_\infty \, \lip(g) \, \theta_k \; .
\end{align*}

\begin{lemma}
  \label{lem:WeakDependenceForCox}
  Let $N$ be a Cox process with stochastic intensity $\xi$, where
  \[
  \xi(A) = \int_A w(s) \rmd s \; ,
  \]
  and $w$ is a positive stationary stochastic process with almost surely
  continuous paths. Assume that $w$ is $\theta$-weak dependent. Then
  \begin{enumerate}[(i)]
  \item \label{item:inutile} the increments $\{\Delta N(k), k\in\Zset\}$ is
    $\theta$-weak dependent with the same rate function as the process $w$;
  \item \label{item:returns:theta-weak-dep} If the efficient shocks $e_i$ are
    i.i.d. and independent of the point process $N$, the sequence of
    calendar-time returns $\{r_k, k \in \mathbb{N}^*\}$, defined by
    \[
    r_k= \sum \limits_{n=N(k-1)+1}^{N(k)} e_n \; ,
    \]
is $\theta$-weak dependent, with the same rate as $w$.
\end{enumerate}
\end{lemma}

\begin{proof} Since~\eqref{item:inutile} is a particular case
  of~\eqref{item:returns:theta-weak-dep} with the shocks taken to be uniformly
  equal to 1, we only need to prove~\eqref{item:returns:theta-weak-dep}.  Let
  $f\in\mcg_u$ and $g\in\mcg_v$.  Consider now the sequence of returns
  $\{r_k,k\in\mathbb{N}^*\}$. By conditioning on the stochastic intensity, we have
\begin{align*}
  \cov(f(r_{i_1}, \dots, r_{i_u}), g(r_{j_1}, \dots, r_{j_v}))
  &= \esp \left[ \cov(f(r_{i_1}, \dots, r_{i_u}), g(r_{j_1}, \dots, r_{j_v})) \mid \xi \right] \\
  &+ \cov \left( \esp(f(r_{i_1}, \dots, r_{i_u})|\xi), \esp( g(r_{j_1}, \dots,
    r_{j_v})|\xi) \right)
\end{align*}
Since the intervals are non intersecting, and $\{e_k\}$ is independent of $N$
which is a Poisson point process conditionally on $\xi$, we have
\[
\esp \left( \cov(f(r_{i_1}, \dots, r_{i_u}), g(r_{j_i}, \dots, r_{j_v})|\xi)
\right)=0 \; .
\]
To compute the second term, denote $h_n(x) = \frac {e^{-x}x^n} {n!}$ and let
$\{e_k^{(i)}\}$, $i\in\mathbb{N}$ be i.i.d. copies of the sequence $\{e_k\}$ and
define $S_i(0)=0$ and $S_i(n) = \sum_{k=1}^n e_k^{(i)}$. Denote also $\xi_i =
\xi(i-1,i]$. Then, by the independent increments property, we have
\begin{align*}
  \esp \left[ f(r_{i_1}, \dots, r_{i_u})|\xi \right] & = \esp \left[ f\left(
      \sum_{k=N(i_1-1)+1}^{N(i_1)} e_k, \dots , \sum_{k=N(i_u-1)+1}^{N(i_u)} e_k
    \right) \middle|\xi \right]
  \\
  & = \sum_{n_1,\dots,n_u=0}^\infty \esp[f(S_{1}(n_1),\dots,S_u(n_u))] \prod_{i=1}^n h_{n_i}(\xi_i) \; .
\end{align*}
Similarly, we have
\[
\esp \left(g(r_{j_1}, \dots, r_{j_v})|\xi \right) = \sum\limits_{m_1, \dots, m_v =
  0}^{\infty} \esp \left[g(S_{1}(m_1), \dots, S_{v}(m_v)) \right] \prod
\limits_{q=1}^{v} h_{m_q}(\xi_{j_q}) \; .
\]
Define the functions $F$ and $G$ on $\Rset^u$ and $\Rset^v$ respectively by
\begin{align*}
  F(x_1,\dots,x_u) & = \sum_{n_1, \dots, n_u = 0}^{\infty} \esp \left[ f(
    S_{1}(n_1), \dots, S_{u}(n_u) ) \right] \prod \limits_{i=1}^{u}  h_{n_i}(x_i)  \; , \\
  G(x_{1}, \dots,x_{v}) & = \sum_{m_1, \dots, m_v = 0}^{\infty} \esp \left[
    g(S_{1}(m_1), \dots, S_{v}(m_v)) \right] \prod_{i=1}^{v} h_{m_i}(x_{i})
\end{align*}
Since $f$ is bounded, we obtain, for all $\boldsymbol{x}\in\Rset^u$,
\[
F(\boldsymbol{x}) \leq \|f\|_\infty \sum_{n_1, \dots, n_u = 0}^{\infty}  \prod_{p=1}^{u} h_{n_p}(x_{p})=||f||_\infty \; .
\]
Thus $\|F\|_\infty \leq \|f\|_\infty$.
We now prove that $F$ is a Lipschitz function. For
$\boldsymbol{x}=(x_1,\dots,x_u)\in\Rset^u$, we have
\begin{align*}
  \frac{\partial F(\boldsymbol{x})}{\partial x_{1}} & = \sum\limits_{n_1=1, n_2,
    \dots, n_u = 0}^{\infty} \esp \left[ f(S_{1}(n_1), \dots, S_{u}(n_u))
  \right] \frac{\rme^{-x_{1}} n_1 x_{1}^{n_1-1}}{n_1!}  \prod_{p=2}^{u}  h_{n_p}(x_{p})  \\
  & \phantom{ = } - \sum_{n_1,\dots, n_u = 0}^{\infty} \esp \left[
    f(S_1(n_1),S_{2}(n_2), \dots , S_{u}(n_u)) \right]  \prod_{p=1}^{u}  h_{n_p}(x_{p})  \\
  & = \sum_{n_1, \dots , n_u = 0}^{\infty} \esp \left[ f(S_{1}(n_1+1),\dots ,
    S_{i_u}(n_u)) \right] \prod_{p=1}^{u} h_{n_p}(x_{p})  \\
  & \phantom{ = } - \sum_{n_1, n_2, \dots , n_u = 0}^{\infty} \esp \left[
    f(S_{1}(n_1), \dots, S_{u}(n_u)) \right] \prod_{p=1}^{u} h_{n_p}(x_{p})
  \\
  & = \sum_{n_1, \dots, n_u = 0}^{\infty} \esp \left[ f(S_{1}(n_1+1),
    S_{2}(n_2), \dots, S_{u}(n_u)) - f(S_{1}(n_1),\dots, S_{u}(n_u)) \right]
  \prod_{p=1}^{u} h_{n_p}(x_{p}) \; .
\end{align*}
Applying the Lipschitz property of $f$, we have
\begin{align*}
  \left| \frac{\partial F(\boldsymbol{x})}{\partial x_{1}} \right| & \leq
  \sum_{n_1, \dots , n_u = 0}^{\infty} \esp \left[ \left| f(S_{1}(n_1+1),
      S_{i_2}(n_2), \dots , S_{i_u}(n_u)) - f(S_{i_1}(n_1), \dots ,
      S_{i_u}(n_u))    \right| \right] \prod_{p=1}^{u} h_{n_p}(x_{p})  \\
  & \leq \lip(f) \sum_{n_1, \dots , n_u = 0}^{\infty} \esp
  \left[|S_{i_1}(n_1+1)-S_{i_1}(n_1)| \right] \prod_{p=1}^{u} h_{n_p}(x_{p})  \\
  & = \lip(f) \esp[|e_{1}|] \sum_{n_1, \dots , n_u = 0}^{\infty} \prod_{p=1}^{u}
  h_{n_p}(x_{p}) = \mathrm{Lip}(f) \esp[|e_1|] \; .
\end{align*}
Similarly, for each $p=2, \dots, u$,
\[
\left \|\frac{\partial F}{\partial x_{p}} \right\|_\infty \leq \mathrm{Lip}(f) \, \esp[|e_1|] \; .
\]
Thus $F$ is also a Lipschitz function and $\mathrm{Lip}(F) \leq \esp[|e_1|]
\mathrm{Lip}(f)$. We can similarly prove that $G$ is a bounded Lipschitz function.  Let us assume for the moment that $\{\xi_i,i \in
\Zset\}$ is $\theta$-weak dependent. We then conclude that
\begin{align*}
  \cov(f(r_{i_1}, \dots , r_{i_u}) , g(r_{j_i}, \dots , r_{j_v})) & = \cov
  \left( \esp[f(r_{i_1}, \dots , r_{i_u})|\xi], \esp[ g(r_{j_i}, \dots,     r_{j_v})\mid \xi] \right) \\
  & = \cov(F(\xi_{i_1}, \dots , \xi_{i_u}), G(\xi_{j_1}, \dots, \xi_{j_v})) \\
  & \leq v \|F\|_\infty \mathrm{Lip}(G) \theta_r \leq v \, \esp[|e_1|] \,
  \|f\|_\infty \, \lip(g) \, \theta_r \; .
\end{align*}
This proves that $\{r_k, k \in \mathbb{N}^*\}$ is $\theta$-weak dependent with the
same rate as $\{\xi_i,i\in\Zset\}$.

Let us now prove the $\theta$-weak dependence of the sequence
$\{\xi_i,i\in\Zset\}$. Let $f\in\mcg_u$ and $g\in\mcg_v$. We must prove that
\begin{align*}
  |\cov(f(\xi_{i_1},\dots,\xi_{i_u}),g(\xi_{j_i},\dots,\xi_{j_v}))| \leq C \, v \,
  \|f\|_\infty \, \lip(g) \theta_r \; ,
\end{align*}
for some constant $C$ and all $(u+v)$-tuples of integers $i_1 < \cdots < i_u <
j_1 < \dots < j_v$ such that $j_1 > i_u+r$.
Since the density $w$ of the stochastic intensity has almost surely continuous
path, it holds that $\xi_i$ is the almost sure limit of the Riemann sums
$\xi_{n,i} = n^{-1}\sum_{q=1}^n w(i-1+q/n)$. Since $f$ and $g$ are continuous and
bounded, we obtain by bounded convergence
\begin{align}
   \label{eq:riemann-sum}
   \cov(f(\xi_{i_1}\dots,\xi_{i_u}),g(\xi_{j_1},\dots,\xi_{j_v})) =
   \lim_{n\to\infty} \cov(f(\xi_{i_1,n}\dots,\xi_{i_u,n}),
   g(\xi_{j_1,n},\dots,\xi_{j_v,n})) \; .
\end{align}
Given a function $h:\Rset^p\to\Rset$, define the function $h_n$ on $\Rset^{np}$ by
\begin{align*}
  h_n(x_{i,q},1 \leq i \leq p,1 \leq q \leq n) = h\left(
    \frac{x_{1,1}+\cdots+x_{1,n}}n,\dots, \frac{x_{p,1}+\cdots+x_{p,n}}n\right)
\end{align*}
If $h$ is bounded and Lipschitz, then $\|h_n\|_\infty \leq \|h\|_\infty$ and
$\lip(h_n) \leq n^{-1} \lip(h)$.  Since the process $w$ is $\theta$ weakly
dependent, we obtain
\begin{align}
  | \cov & (f(\xi_{i_1,n}\dots,\xi_{i_u,n}), g(\xi_{j_1,n},\dots,\xi_{j_v,n})) | \nonumber \\
  & = |\cov(f_n(w(i_s-1+\frac qn),s=1,\dots,u,q=1,\dots,n), g_n(w(j_t-1+\frac qn),t=1,\dots,v,q=1,\dots,n)) | \nonumber \\
  & \leq vn \times \|f\|_\infty \times  n^{-1} \lip(g) \times \theta_{k-1} \nonumber  \\
  & = v \|f\|_\infty \, \lip(g) \theta_{k-1} \; .
 \label{eq:theta-dep-w}
\end{align}
The bounds~(\ref{eq:riemann-sum}) and~(\ref{eq:theta-dep-w}) proves that the
sequence $\{\xi_i,i\in\Zset\}$ is $\theta$ weak dependent with rate $\theta_{k-1}$.
\end{proof}

\end{document}